\documentclass[11pt]{article}
\usepackage{amsmath,amsfonts,amssymb}
\usepackage{latexsym}
\usepackage{epsfig,overpic}
\usepackage{subfigure}
\usepackage{color}

\newcommand{\bn}{\mathbf n}
\newcommand{\bp}{\mathbf p}

\newcommand{\bw}{\mathbf w}
\newcommand{\bx}{\mathbf x}
\newcommand{\by}{\mathbf y}

\newcommand{\F}{\mathcal F}

\newcommand{\T}{\mathcal T}

\newcommand{\Div}{\mathop{\rm div}}

\newcommand{\Gs}{\mathcal{G}} 
\DeclareGraphicsExtensions{.pdf,.eps,.ps,.eps.gz,.ps.gz,.eps.Y}

\def\ds{\mathrm{d}\mathbf{s}}

\def\nath{\nabla_{\Gamma_h}}

\begin{document}
\title{A trace finite element method for PDEs on evolving surfaces
\thanks{Partially supported by NSF through the Division of Mathematical Sciences grant  1522252.}}
\author{Maxim A. Olshanskii\thanks{Department of Mathematics, University of Houston, Houston, Texas 77204-3008
(molshan@math.uh.edu)}
\and Xianmin Xu\thanks{LSEC, Institute of Computational Mathematics and Scientific/Engineering Computing,
  NCMIS, AMSS, Chinese Academy of Sciences, Beijing 100190, China (xmxu@lsec.cc.ac.cn).}
}
\date{}

\maketitle
\begin{abstract} In this paper, we propose an approach for solving PDEs on evolving surfaces using a combination of  the trace finite element method and a fast marching method.
The numerical approach is based on the Eulerian description of the surface problem and employs a time-independent background mesh that is not fitted to the surface. The surface and its evolution may be given implicitly, for example, by the level set method. Extension of the PDE off the surface is \textit{not} required. The method introduced in this paper naturally allows a surface to undergo topological changes and experience local geometric singularities.
In the simplest setting, the numerical method is second order accurate in space and time. Higher order variants are feasible, but not studied in this paper. We show results of several numerical experiments, which demonstrate the convergence properties of the method and its ability to handle the case of the surface with topological changes.
\end{abstract}
%
\section{Introduction}\label{sec:introduction}
Partial differential equations on evolving surfaces arise in a number of  mathematical models in natural sciences and engineering.   Well-known examples include the diffusion and transport of  surfactants along interfaces in multiphase fluids \cite{GReusken2011,Milliken,Stone}, diffusion-induced grain boundary motion \cite{GrainBnd1,GrainBnd2} and lipid interactions in moving cell membranes \cite{ElliotStinner,Novaketal}. Thus, recently there has been a significant  interest in developing and analyzing numerical methods for PDEs on time-dependent surfaces. While all of finite difference, finite volumes and finite element methods have been considered in the literature for numerical solution of PDEs on manifolds, in this work we focus on finite element methods.

The choice of a numerical approach for solving a PDE on evolving surface $\Gamma(t)$ largely depends on which of Lagrangian or Euclidian frameworks is used to setup the problem and describe the surface evolution.
In \cite{Dziuk07,DziukElliot2013a,elliott2015error} Elliott and co-workers developed and analyzed a finite element method (FEM) for computing transport and diffusion on a surface which is based on a Lagrangian tracking of the surface evolution. Some recent developments of the finite element methods for surface PDEs   based on the Lagrangian description can be found, e.g., in \cite{barrett2015stable,bretschneider2016solving,elliott2015evolving,macdonald2016computational,sokolov2015afc}.
If a surface undergoes strong deformations, topological changes, or it is defined implicitly, e.g., as the zero level of a level set function, then numerical methods based on the Lagrangian approach have certain disadvantages. Methods using an Eulerian approach were developed in \cite{AS03,bertalmio2001variational,xu2006level,XuZh}, based on an extension  of the surface PDE into a bulk domain that contains the surface. Although in the original papers, finite differences were used, the approach is also suitable for finite element methods, see, e.g., \cite{burger2009finite} 
Related technique is the closest point method in \cite{petras2016pdes}, where the closest point representation of the surface and differential operators is used in an ambient space to allow a standard Cartesian finite difference discretization method. 

In the present paper, we develop yet another finite element method for solving a PDE on a time-dependent surface $\Gamma(t)$. The surface is embedded in a bulk computational domain. We assume a sharp representation of the surface rather than a diffusive interface approach typical for the phase-field  models of interfacial problems. The level set method \cite{SethianBook} is suitable for the purposes of this paper and will be used here to recover the evolution of the surface. 
We are interested in a surface FEM known in the literature as the trace or cut FEM.  The trace finite element  method uses the restrictions (traces) of a function from  the background time-independent finite element space on the reconstructed discrete surface. This does not involve any mesh fitting towards the surface or an extension of the PDE.

The trace FEM method was originally introduced for elliptic PDEs on stationary surfaces in \cite{OlshReusken08}. Further the analysis and several extensions of the method were developed in the series of publications. This includes higher order,  stabilized, discontinuous Galerkin and adaptive variants of the method as well as applications to the surface--bulk coupled transport--diffusion problem, two-phase fluids with soluble surfactants and coupled bulk-membrane elasticity problems, see, e.g., \cite{Alg1,burman2016cutb,burman2016full,cenanovic2015cut,chernyshenko2015adaptive,DemlowOlsh,grande2016higher,
gross2015trace,lehrenfeld2016high,OlsR2009,ORXimanum,reusken2015analysis}. There have been several successful attempts to extend the  method to time-dependent surfaces.
In \cite{deckelnick2014unfitted} the trace FEM was combined with the narrow-band unfitted FEM from \cite{deckelnick2009h} to devise an unfitted finite element method for parabolic equations on evolving surfaces. The resulting method preserves mass in the case of an advection-diffusion conservation law. The method based on a characteristic-Galerkin formulation combined with the trace
FEM in space was proposed in \cite{hansbo2015characteristic}. Thanks to the semi-Lagrangian treatment of the material derivative {\color{black}(numerical integration back in time along characteristics)} this variant of the method does not require stabilization for the dominating advection. The first order convergence of the characteristic--trace FEM was demonstrated by a rigorous \textit{a priori} error analysis and in numerical experiments. Another direction was taken in  \cite{olshanskii2014eulerian}, where a space--time weak formulation of the surface problem was introduced. Based on this weak formulation, space--time variants of the trace FEM for PDEs on evolving surfaces were proposed in that paper and in \cite{grande2014eulerian}.
The method from \cite{olshanskii2014eulerian} employs  discontinuous piecewise linear in time -- continuous piecewise linear in space finite elements.   In \cite{olshanskii2014error} the first order convergence in space and time of the method in an energy norm and second order  convergence in a weaker norm was proved. In \cite{grande2014eulerian}, the author experimented with both  continuous and discontinuous  in time piecewise linear finite elements.

In the space--time trace FEM, the  trial and test finite element spaces consist of traces of standard volumetric elements on a space--time manifold resulting from the evolution of a surface. The implementation requires the numerical integration over the tetrahedral reconstruction of
the 3D manifold embedded in the  $\mathbb{R}^4$ ambient space. An efficient algorithm for such numerical reconstruction was suggested in \cite{grande2014eulerian} and implemented in the DROPS finite element package \cite{DROPS}. In \cite{hansbo2016cut} a stabilized version of the space--time trace FEM for coupled bulk-surface problems was implemented using Gauss-Lobatto quadrature rules in time. In this implementation, one does not reconstruct the 3D space--time manifold but instead needs the 2D surface approximations in the quadrature nodes. The numerical experience with space--time trace FEM based on quadrature rules in time is mixed. The authors of \cite{hansbo2016cut} reported a second order convergence of the method for a number of 2D tests (in this case a 1D PDE is posed on an evolving curve), while in \cite{grande2014eulerian} one finds an example of a  2D smoothly  deforming surface  when the space--time method based on the trapezoidal  quadrature rule fails to deliver convergent results. The error analysis of such simplified versions is an open question.

Although the space--time framework is natural for the development of unfitted FEMs for PDEs on evolving surfaces,
the implementation of such methods is not straightforward, especially if a higher order method is desired.
In this paper, we propose a variant of the trace FEM for time-dependent surfaces that uses simple finite difference approximations
of \textit{time} derivatives. It avoids any reconstruction of the surface--time manifold, it also avoids finding surface approximations at quadrature nodes. Instead, the method requires arbitrary, but smooth in a sense clarified later, extension of the numerical solution
off the surface to a narrow strip around the surface. We stress that in the present method one does not extend either problem data or differential operators to a surface neighborhood as in the methods based on PDEs extension.
At a given time node $t_n$, the degrees of freedom in the narrow strip (except those belonging to tetrahedra cut by the surface $\Gamma(t_n)$) do not contribute to algebraic systems, but are only used to store the solution values from several previous time steps. In numerical
examples, we use the BDF2 scheme for time integration and so the narrow band degrees of freedom store the finite element solution for $t=t_{n-1}$ and $t=t_{n-2}$. To find a suitable extension, we apply a variant of the Fast Marching Method (FMM), see, e.g., \cite{GReusken2011,sethian1996fast}. At each time step, the trace FEM for a PDE on a steady surface $\Gamma(t_n)$ and the FMM are used in a modular way, which makes the implementation straightforward in a standard or legacy finite element software. For P1 background finite elements and BDF2 time stepping scheme, {\color{black} numerical experiments show that }
the method is second order accurate (assuming $\Delta t\simeq h$) and has no stability restrictions on the time step. We remark that the numerical methodology naturally extends to the surface-bulk coupled problems with propagating
interfaces. However, in this paper we concentrate on the case  when surface processes are decoupled from  processes in the bulk.

The remainder of the paper is organized as follows. In section \ref{s_setup} we present the PDE model on an evolving surface
and review some properties of the model. Section~\ref{s_FEM} introduces our variant of the trace FEM, which avoids space--time elements. Here we discuss implementation details. Section~\ref{s_num} collects the results for a series of numerical experiments. The experiments aim to access the accuracy of the method as well as the ability to solve PDEs along a surface undergoing topological changes. For the latter purpose we consider the example of  the diffusion of a surfactant on a surface of two colliding droplets.

\section{Mathematical formulation}\label{s_setup}

Consider a surface $\Gamma(t)$ passively advected by a smooth velocity field $\bw=\bw(\bx,t)$, i.e. the normal velocity of $\Gamma(t)$ is given by $\bw \cdot \bn$, with
$\bn$ the unit normal on $\Gamma(t)$. We assume that for all $t \in [0,T] $,  $\Gamma(t)$ is a smooth hypersurface that is  closed ($\partial \Gamma =\emptyset$), connected, oriented, and contained in a fixed domain $\Omega \subset \Bbb{R}^d$, $d=2,3$. In the remainder we consider $d=3$, but all results have analogs for the case $d=2$.

As an example of the surface  PDE, consider the transport--diffusion equation modelling the conservation of a scalar quantity
$u$ with a diffusive flux on $\Gamma(t)$ (cf. \cite{James04}):
\begin{equation}
\dot{u} + ({\Div}_\Gamma\bw)u -{ \nu}\Delta_{\Gamma} u=0\quad\text{on}~~\Gamma(t), ~~t\in (0,T],
\label{transport}
\end{equation}
 with initial condition $u(\bx,0)=u_0(\bx)$ for $\bx \in \Gamma_0:=\Gamma(0)$.
 Here $\dot{u}$ denotes the advective material derivative, ${\Div}_\Gamma:=\operatorname{tr}\left( (I-\bn\bn^T)\nabla\right)$ is the surface divergence,  $\Delta_\Gamma$ is the  Laplace--Beltrami operator, and $\nu>0$ is the constant diffusion coefficient. The well-posedness of suitable weak formulations of \eqref{transport}  has been proved in {\color{black}\cite{Dziuk07,olshanskii2014eulerian,alphonse2014abstract}.}

The equation \eqref{transport} can be written in several equivalent forms, see \cite{DEreview}.  In particular,
for any smooth extension of   $u$ from the space--time manifold
\[
\Gs: = \bigcup\limits_{t \in (0,T)} \Gamma(t) \times \{t\},\quad  \Gs\subset \Bbb{R}^{4},
\]
 to a neighborhood of $\Gs$, one can expand the material derivative $\dot{u}= \frac{\partial u}{\partial t} + \bw \cdot \nabla u$.
 Note that the identity holds independently of a smooth extension of $u$ off the surface.

Assume further that the surface is defined implicitly as the zero level of the smooth 
level set function $\phi$ on $\Omega\times(0,T)$:
\[
\Gamma(t)=\{\bx\in\mathbb{R}^3\,:\,\phi(t,\bx)=0\},
\]
such that $|\nabla\phi|\ge c_0>0$ in $\mathcal{O}(\Gs)$, a  neighborhood of $\Gs$. One can consider an extension  $u^e$ in  $\mathcal{O}(\Gs)$ such that $u^e=u$ on $\Gs$ and $\nabla u^e\cdot\nabla \phi =0$ in $\mathcal{O}(\Gs)$.
Note that $u^e$ is smooth once $\phi$  and $u$ are both smooth. Below we use the same notation
$u$ for the solution of the surface PDE  \eqref{transport} and its extension. We have the equivalent formulation,
\begin{equation}
\left\{\begin{split}
 \frac{\partial u}{\partial t} + \bw \cdot \nabla u + ({\Div}_\Gamma\bw)u -{ \nu}\Delta_{\Gamma} u&=0\qquad\text{on}~~\Gamma(t), \\
 \nabla u\cdot\nabla \phi& =0 \qquad\text{in}~\mathcal{O}(\Gamma(t)),
 \end{split}
 \right.~~t\in (0,T].
\label{transport_new}
\end{equation}
If $\phi$ is the signed distance function, the second equation in \eqref{transport_new} defines the normal extension of $u$, i.e. the solution $u$ stays constant in normal directions to $\Gamma(t)$.
{\color{black}Otherwise, $\nabla u\cdot\nabla \phi=0$ defines an extension, which is not necessarily the normal extension. In fact, any  extension is suitable for our purposes, if $u$ is smooth function in a neighborhood of $\Gs$.  We shall make an exception is section~\ref{s_tFEMs}, where error analysis is reviewed and  we need the normal extension to formulate certain estimates.}

In the next section, we devise the trace FEM based on  the formulation  \eqref{transport_new}.

\section{The finite element method}\label{s_FEM} We first collect some preliminaries and recall the trace FEM from \cite{OlshReusken08} for the elliptic equations on stationary surfaces and some of its properties. Further, in section~\ref{s_tFEMe} we apply this method on each time step of a numerical algorithm for the transient problem \eqref{transport_new}.

\subsection{Background mesh and induced surface triangulations}\label{s_surf}
Consider a tetrahedral subdivision $\mathcal{T}_h$ of the bulk computational domain $\Omega$. We assume that the triangulation
$\mathcal{T}_h$ is regular (no hanging nodes). For each tetrahedron $S\in \mathcal{T}_h$, let $h_S$ denote its diameter and define the global parameter of the triangulation by $h = \max_{S} h_S$. We assume that $\mathcal{T}_h$ is shape regular, i.e. there exists $\kappa>0$ such that for every $S\in \mathcal{T}_h$  the radius $\rho_S$ of its inscribed sphere satisfies
\begin{equation}\label{shaperegularity}
\rho_S>h_S/\kappa.
\end{equation}

For each time $t\in[0,T]$, denote by $\Gamma_h(t)$ a polygonal approximation of $\Gamma(t)$. We assume that
$\Gamma_h(t)$ is a $C^{0,1}$ surface without  boundary and $\Gamma_h(t)$ can be partitioned in planar triangular segments:
\begin{equation} \label{defgammah}
 \Gamma_h(t)=\bigcup\limits_{T\in\mathcal{F}_h(t)} T,
\end{equation}
where $\mathcal{F}_h(t)$ is the set of all  triangular segments $T$.
We assume that for any $T\in\mathcal{F}_h(t)$ there is only \textit{one} tetrahedron $S_T\in\mathcal{T}_h$ such that $T\subset S_T$
(if $T$ lies on a face shared by two tetrahedra, any of these two tetrahedra can be chosen as $S_T$).

For the level set description of $\Gamma(t)$, the polygonal surface $\Gamma_h(t)$ is defined by the finite element level set function as follows. Consider a continuous function $\phi_h(t,\bx)$ such that for any $t\in[0,T]$ the function $\phi_h$ is piecewise linear with respect
to the triangulation $\mathcal{T}_h$.  Its zero level set defines $\Gamma_h(t)$,
\begin{equation}\label{Gammah}
\Gamma_h(t):=\{\bx\in\Omega\,:\, \phi_h(t,\bx)=0 \}.
\end{equation}
We assume that $\Gamma_h(t)$ is an approximation to $\Gamma(t)$. This  is  a reasonable assumption if $\phi_h$ is
either an interpolant to the known $\phi$ or one finds $\phi_h$ as the solution to a discrete level set equation. In the later case, one may have no direct knowledge of $\Gamma(t)$. Other interface capturing techniques such as the volume of fluid method also can be used subject to a postprocessing step to recover $\Gamma_h$.

 \begin{figure}[ht!]
 \centering
  \subfigure[]{
   \includegraphics[width=2.4in]{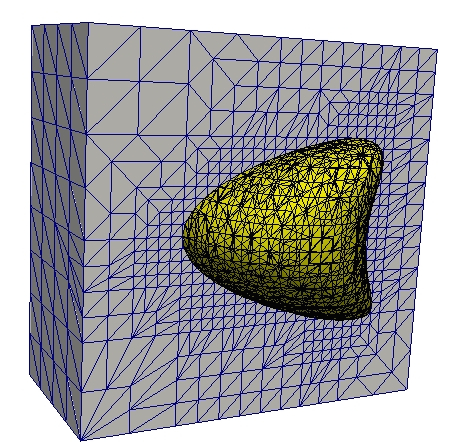}
   }
   \hspace{0.6cm}
  \subfigure[]{
  \includegraphics[width=2.5in]{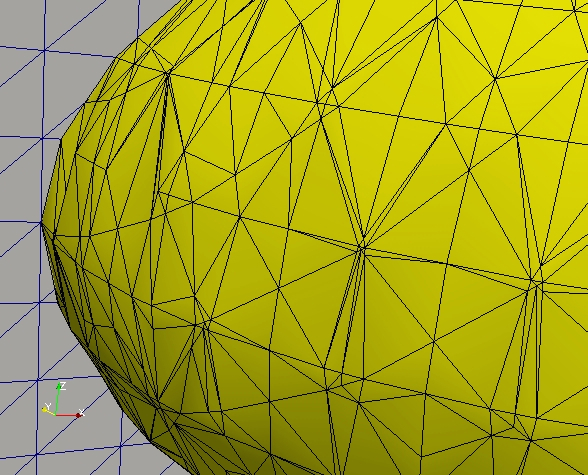}
}
 \caption{Left: Cut of the background   and induced surface meshes for $\Gamma_h(0)$ from Experiment~4 in Section~\ref{s_num}. Right: The zoom-in of the surface mesh.}
 \label{fig:surfaceMesh}
\end{figure}

The intersection of $\Gamma_h(t)$ defined in \eqref{Gammah} with any tetrahedron in $\T_h$ is either a triangle or a quadrilateral.
{\color{black} If the intersection is a quadrilateral, we divide it into two triangles.} This construction of $\Gamma_h(t)$ satisfies the assumptions made above. The bulk triangulation $\T_h$ consisting of tetrahedra and the induced surface triangulation  are illustrated in Figure~\ref{fig:surfaceMesh}.
There are no restrictions on how $\Gamma_h(t)$ cuts through the background mesh, and so for any fixed time instance $t$ the resulting triangulation $\mathcal{F}_h(t)$  is \emph{not} necessarily  regular. The elements from $\mathcal{F}_h(t)$ may have very small internal angles and the size of neighboring triangles can vary strongly, cf.~Figure~\ref{fig:surfaceMesh} (right).
Thus, {\color{black} $\Gamma_h(t)$ is not a regular triangulation of $\Gamma(t)$.}
Two interesting properties of the induced surface triangulations are known in the literature ~\cite{DemlowOlsh,olshanskii2012surface}: (i)~If the background triangulation $\T_h$ satisfies the minimal angle condition~\eqref{shaperegularity},  then the surface triangulation satisfies \textit{maximum} angle condition; (ii)~Any element from  $\mathcal{F}_h(t)$ shares at least one vertex with a full size shape regular triangle from $\mathcal{F}_h(t)$. The trace finite element method does not exploit these properties directly, but they are still useful if one is interested in understanding the performance of the method.

 We note that the surface triangulations $\mathcal{F}_h(t)$ will be used only to perform numerical integration in the finite element method below, while approximation properties of the method, as we shall see, depend on the volumetric tetrahedral mesh.

\subsection{The trace FEM: steady surface}\label{s_tFEMs} To review the idea of the trace FEM, assume for a moment the stationary transport--diffusion problem on a steady  closed smooth surface $\Gamma$,
\begin{equation}
\alpha u+\bw \cdot \nabla u + ({\Div}_\Gamma\bw)u -{ \nu}\Delta_{\Gamma} u=f\quad\text{on}~~\Gamma.
\label{transport_st}
\end{equation}
Here we assume $\alpha>0$ and $\bw\cdot\bn=0$. Integration by parts over $\Gamma$ gives the weak formulation of \eqref{transport_st}:
Find $u\in H^1(\Gamma)$ such that
\begin{equation}
\int_{\Gamma}\left(\alpha \, uv+\nu\nabla_\Gamma u\cdot\nabla_\Gamma v\, - (\bw\cdot\nabla  v) u\,\right)\, \ds =\int_{\Gamma}f v\, \ds \label{weak}
\end{equation}
for all $v\in H^1(\Gamma)$.  In the trace FEM, one substitutes $\Gamma$ with $\Gamma_h$ in \eqref{weak} {\color{black}($\Gamma_h$ is constructed as in section~\ref{s_surf})}, and instead of $H^1(\Gamma)$ considers the space of traces on $\Gamma_h$ of all functions from the background ambient finite element space.
This can be formally defined as follows.

Consider the volumetric finite element space of all piecewise linear continuous functions with respect to the bulk triangulation $\mathcal{T}_h$:
\begin{equation}
 V_h:=\{v_h\in C(\Omega)\ |\ v|_{S}\in P_1~~ \forall\ S \in\mathcal{T}_h\}.
 \label{e:2.6}
\end{equation}
$V_h$ induces the following space on $\Gamma_h$:
\begin{equation*}
 V_h^{\Gamma}:=\{\psi_h\in C(\Gamma_h)\ |\ \exists ~ v_h\in V_h\  \text{such that }\ \psi_h=v_h~\text{on}~{\Gamma_h}\}.
\end{equation*}

Given the surface finite element space $V_h^{\Gamma}$,  the finite element
discretization of \eqref{transport_st} reads:   Find $u_h\in V_h^{\Gamma}$ such that
\begin{equation}
\int_{\Gamma_h}\left(\alpha \, u_hv_h+\nu\nath u_h\cdot\nath v_h\, - (\bw\cdot\nabla v_h) u_h
\right)\, \ds_h =\int_{\Gamma_h}f_h v_h\, \ds_h \label{FEM}
\end{equation}
for all $v_h\in V_h^{\Gamma}$. Here $f_h$ is an approximation of the problem source term on $\Gamma_h$.

 {\color{black}From here and up to the end of this section} $f^e$ denotes a \emph{normal} extension of a quantity $f$ from $\Gamma$. For a smooth closed surface, $f^e$ is well defined in a neighborhood $\mathcal{O}(\Gamma)$.
Assume that $\Gamma_h$ approximates $\Gamma$ in the following sense: It holds $\Gamma_h\subset\mathcal{O}(\Gamma)$ and
\begin{equation}
\|\bx-\bp(\bx)\|_{L^\infty(\Gamma_h)}+h\|\bn^e-\bn_h\|_{L^\infty(\Gamma_h)}\le c\,h^2,
\label{G_Gh}
\end{equation}
where $\bn_h$ is an external normal vector on $\Gamma_h$ and $\bp(\bx)\in\Gamma$ is the closest surface point for $\bx$. Given \eqref{G_Gh}, the trace FEM  is second order accurate in the $L^2$ surface norm and first order accurate in $H^1$ surface norm~\cite{OlshReusken08,ORXimanum}: For solutions of \eqref{transport_st} and \eqref{FEM}, it holds
\[
\|u^e-u_h\|_{L^2(\Gamma_h)}+h\|\nabla_{\Gamma_h}(u^e-u_h)\|_{L^2(\Gamma_h)} \le c\,h^2,
\]
with a constant $c$ dependent only on the shape regularity of $\T_h$  and \textit{independent of how the
surface $\Gamma_h$ cuts through the background mesh}. This robustness property is extremely useful for extending the method to time-dependent  surfaces. It allows to keep the same background mesh while the surface evolves through the bulk domain, avoiding unnecessary mesh fitting and mesh reconstruction.

Before we consider the time-dependent case, a few  important properties of the method should be mentioned.
Firstly, the authors of \cite{deckelnick2014unfitted} noted that one can use the full gradient instead of the  tangential gradient in the diffusion term in \eqref{FEM}. This leads to the following FEM formulation:
 Find $u_h\in V_h$ such that
\begin{equation}
\int_{\Gamma_h}\left(\alpha \, u_hv_h+\nu\nabla u_h\cdot\nabla v_h\, - (\bw\cdot\nabla v_h) u_h
\right)\, \ds_h =\int_{\Gamma_h}f_h v_h\, \ds_h \label{FEMh}
\end{equation}
for all $v_h\in V_h$. The rationality behind the modification is clear from the following observation.
For the normal extension $u^e$ of the solution $u$  we have $\nabla_{\Gamma} u= \nabla u^e$ and so $u^e$ satisfies the integral equality \eqref{weak} with surface gradients (in the diffusion term) replaced
by full gradients and for arbitrary smooth function $v$ on $\Omega$. 
Therefore, by solving  \eqref{FEMh} we recover  $u_h$, which approximates
the PDE solution $u$ on the triangulated surface $\Gamma_h$ \textit{and} its normal extension $u^e$ in the strip consisting of all tetrahedra cut by the surface $\Gamma_h$.

The formulation  \eqref{FEMh} uses the bulk finite element space $V_h$ instead of the surface finite element space $V_h^\Gamma$ in \eqref{FEM}. However, the practical implementation of both methods uses the same frame of all bulk finite element nodal basis functions $\phi_i\in V_h$ such that  $\mbox{supp}(\phi_i)\cap\Gamma_h\neq\emptyset$. Hence, the active degrees of freedom  in both methods are the same. The stiffness matrices are, however, different.
For the case of the Laplace-Beltrami problem and a  regular quasi-uniform tetrahedral grid, the studies in \cite{deckelnick2014unfitted,reusken2015analysis}  show that the conditioning of the (diagonally scaled) stiffness matrix of the method \eqref{FEMh}  improves over the conditioning of the matrix for \eqref{FEM}, at the expense of a slight deterioration of the accuracy. Further in this paper we shall use the full gradient version of the trace FEM.

From the formulations  \eqref{FEMh} or \eqref{FEM} we see that only those degrees of freedom of the background finite element
space $V_h$ are active (enter the system of algebraic equations) that are tailored to the tetrahedra cut by $\Gamma_h$.
This provides us with a method of optimal computational complexity, which is not always the case for the methods based on an extension of surface PDE to the bulk domain.  Due to the possible small cuts of bulk tetrahedra (cf. Figure~\ref{fig:surfaceMesh}), the resulting stiffness matrices can be poor conditioned. The simple diagonal re-scaling of the matrices significantly improves the conditioning and eliminates outliers in the spectrum, see \cite{OlsR2009,reusken2015analysis}. Therefore, Krylov subspace iterative
methods applied to the re-scaled matrices are very efficient to solve the algebraic systems.
Since the resulting matrices are sparse and resemble discretizations of 2D PDEs, using an optimized direct solver is also a suitable option. 

\subsection{The trace FEM: evolving surface}\label{s_tFEMe}
For the evolving surface case, we extend the approach in such a way that the trace FEM \eqref{FEMh} is applied on each time
step for the recovered surface $\Gamma_h^{n}\approx\Gamma(t_n)$. Here and further, $\{t_n\}$, with $0=t_0<\dots<t_n<\dots<t_N=T$, is the temporal mesh, and $u^{n}$ approximates $u(t_{n})$. As before, $V_h$ is a \textit{time independent} bulk finite element space with respect to the given background triangulation $\T_h$.

Assume that  a smooth extension $u^e(\bx,t)$  is available in $\mathcal{O}(\Gs)$, and
\begin{equation}\label{Cond1}
\Gamma(t_n)\subset \left\{\bx\in\Omega\,:\,(\bx,t_{n-1})\in\mathcal{O}(\Gs)\right\}.
\end{equation}
In this case, one may discretize \eqref{transport_new} in
time using, for example, the implicit Euler method:
\begin{equation}
\frac{u^{n}-u^e(t_{n-1})}{\Delta t} + \bw^{n} \cdot \nabla u^{n} + ({\Div}_\Gamma\bw^{n})u^{n} -{ \nu}\Delta_{\Gamma} u^{n}=0\quad\text{on}~\Gamma(t_n),
\label{transportFD}
\end{equation}
$\Delta t=t_n-t_{n-1}$.
Now we apply the trace FEM to solve \eqref{transportFD} numerically. The trace FEM is a natural choice here, since $\Gamma(t_n)$ is not fitted by the mesh. We look for $u^{n}_h\in V_h$ solving
\begin{equation}\int_{\Gamma_h^{n}}\left(\frac{1}{\Delta t}u^{n}_hv_h - (\bw^{n}_h \cdot \nabla v_h) u^{n}_h \right)\,\ds_h +{ \nu}\int_{\Gamma_h^{n}}\nabla u^{n}_h\cdot\nabla v_h\,\ds_h =\int_{\Gamma_h^{n}}\frac{1}{\Delta t}u^{e,n-1}_hv_h\,\ds_h
\label{transportFDFE}
\end{equation}
for all $v_h\in V_h$. Here $u^{e,n-1}_h$ is a suitable extension of $u^{n-1}_h$ from $\Gamma_h^{n-1}$ to the surface neighborhood, $\mathcal{O}(\Gamma_h^{n-1})$. Condition \eqref{Cond1} yields to the condition
\begin{equation}\label{Cond2}
\Gamma^{n}_h\subset \mathcal{O}(\Gamma^{n-1}_h).
\end{equation}
Note that \eqref{Cond2} is not a Courant condition on $\Delta t$, but rather a condition on a width of a strip surrounding
the surface, where the extension of the finite element solution is performed.
{\color{black}Over one time step, a material point on the surface can travel a distance not exceeding $\|\mathbf{w}\|_{L^\infty}\Delta t$. Therefore, it is safe to
extend the solution to all tetrahedra intersecting the strip of the width $2\|\mathbf{w}\|_{L^\infty}\Delta t$ surrounding the surface.  Hence, we consider all tetrahedra having at least one vertex closer than $\|\mathbf{w}\|_{L^\infty}\Delta t$ to the surface: Define
\begin{equation}\label{strip}
\mathcal{\widetilde{S}}(\Gamma_h^{n}):=\left\{ S\in \mathcal{T}_h~:~\exists~x\in \mathcal{N}(S),~~\text{s.t.}~~\hbox{dist}(\bx,\Gamma_h^n)<L\|\mathbf{w}\|_{L^\infty}\Delta t\right\},\quad L=1,
\end{equation}
where $\mathcal{N}(S)$ is the set of all nodes for $S\in \mathcal{T}_h$.
The criterion in \eqref{strip} can be refined by exploiting the local information about $\bw$ or about $\bn\cdot\bw$.}

After we determine the numerical extension procedure, $u^{k}_h\to u^{e,k}_h$,  the identity \eqref{transportFDFE} defines the fully discrete numerical method.

In general, to find a suitable extension, one can consider a numerical solver for hyperbolic systems and apply it to the second equation in \eqref{transport_new}.
For example, one can use a finite element method to solve the problem
\[
\frac{\partial u^e}{\partial t'}+  \nabla u^e\cdot\nabla \phi(t^k) =0,\quad\text{such that}~~u^e=u^{k}_h~~\text{on}~\Gamma_h^k,
\]
with the auxiliary time $t'$, and let $u^{k,e}_h:=\lim\limits_{t'\to\infty}u^{e}(t')$.  Another technique to compute extensions
(also used for the re-initialization of the signed distance function in the level-set method) is the Fast Marching Method \cite{sethian1996fast}. We find the FMM technique convenient and fast for building suitable extensions in narrow bands of tetrahedra containing $\Gamma_h$.  We give the details of the FMM  in the next section.

We need one further notation. Denote  by $\mathcal{S}(\Gamma_h^{k})$ a strip of all tetrahedra cut by $\Gamma_h^{k}$:
\[
\mathcal{S}(\Gamma_h^{k})=\bigcup_{S\in\mathcal{T}_{\Gamma}^k} \overline{S},\quad\text{with}~~ \mathcal{T}_{\Gamma}^k := \{S\in \mathcal{T}_h: S \cap \Gamma^k_h \neq \emptyset  \}.
\]
We want to exploit the fact that the trace finite element method  provides  us with the normal extension in $\mathcal{S}(\Gamma_h^{k})$ `for free', since the solution $u^{n}_h$ of \eqref{transportFDFE}  approximately satisfies  $\frac{\partial u_h^k}{\partial \bn}=0$ in $\mathcal{S}(\Gamma_h^{k})$, by the property of the full gradient FEM formulation.

\medskip

For given $u^{e,n-1}_h$ and $\phi_h(t_n)$, one time step of the algorithm now  reads:
\begin{enumerate}
\item Solve \eqref{transportFDFE} for $u^{n}_h\in V_h$;\\[-3ex]
\item Apply the FMM to find $u^{e,n}_h$ in $\mathcal{\widetilde{S}}(\Gamma_h^{n})\setminus \mathcal{S}(\Gamma_h^{n})$
such that
$
 u^{e,n}_h=u^{n}_h~\text{on}~\partial \mathcal{S}(\Gamma_h^{n}).
$
\end{enumerate}
\smallskip
If the motion of the surface is coupled to the solution of the surface PDE (the examples include two-phase flows with surfactant or some models of tumor growth~\cite{GReusken2011,chaplain2001spatio}), then a method to find an evolution of $\phi_h$ has to be added, while finding  $u^{e,n}_h$ can be combined with a re-initialization of $\phi_h$ in the FMM.

A particular advantage of the present variant of the trace FEM for evolving domains is that the accuracy order
in time can be easily increased using standard finite differences. In numerical experiments we use
 the BDF2 scheme: The first term in \eqref{transportFD} is replaced by
 \[
\frac{3u^{n}-4u^e(t_{n-1})+u^e(t_{n-2})}{2\Delta t}
\]
and we set $L=2$ in \eqref{strip}; the corresponding modifications in \eqref{transportFDFE} are obvious.
Furthermore, one may increase the accuracy order  in space by using higher order background finite elements and a higher order surface reconstruction, see \cite{grande2016higher,lehrenfeld2016high} for practical higher order variants of trace FEM on stationary surfaces. In the framework of this paper, the use of these higher order methods is decoupled  from the numerical integration in time.

\subsection{Extension by FMM}\label{s_FMM}
The Fast Marching Method is a well-known technique to  compute an approximate distance function to an interface embedded in a computational domain.
Here we build on the variant of the FMM from section~7.4.1 of \cite{GReusken2011} to compute  finite element function extensions in a strip of tetrahedra. We need some further notations. For a vertex $\bx$ of the background triangulation $\T_h$, $\mathcal{S}(\bx)$ denotes the union of all tetrahedra sharing $\bx$. We fix $t_n$ and
let $\Gamma_h=\Gamma_h^{n}$, $\mathcal{S}(\Gamma_h)=\mathcal{S}(\Gamma_h^{n})$.
Note that we do not necessarily have  \textit{a priori} information of  $\mathcal{\widetilde{S}}(\Gamma_h)$, since the distance
function may not be available. Finding the narrow band for the extension is a part of the FMM below.
We need the set of vertices from tetrahedra cut by the mesh:
\[
\mathcal{N}_{\Gamma}=\{\bx\in\mathbb{R}^3\,:\, \bx\in \mathcal{N}(S)~\hbox{ for some } S\in \mathcal{S}(\Gamma_h) \}.
\]

Assume $u_h=u^{n}_h\in V_h$ solves \eqref{transportFDFE} and we are interested in computing $u_h^e$ in $\mathcal{\widetilde{S}}(\Gamma_h)$.
The FMM is based on a greedy grid traversal technique and  consists of two phases.
\smallskip

{\it Initialization phase.}
In the tetrahedra cut by $\Gamma_h$ the full-gradient trace FEM provides us with the normal extension. Hence, we set
$$
u_h^{e}(\bx)=u_h(\bx)\quad\text{for}~\bx\in\mathcal{N}_{\Gamma}.
$$
For the next step of FMM, we also need a distance function $d(\bx)$ for all $x\in\mathcal{N}_{\Gamma}$.
For any $S_T\in \mathcal{S}(\Gamma_h)$, we know that $T=S_T\cap\Gamma_h$
is a triangle or quadrilateral with vertices  $\{\by_j\}$, $j=1,\dots,J$, where $J=3$ or $J=4$. Denote by $\mathbb{P}_T$ the plane containing $T$ and by $P_h \bx$  the projection of $\bx$ on $\mathbb{P}_T$.
Then, for each $\bx\in \mathcal{N}(T)$, we define
\begin{equation}
d_T(\bx):=\left\{
\begin{array}{ll}
|\bx-P_h \bx|& \hbox{if } P_h \bx\in T,\\
\min\limits_{1\leq j\leq J}|\bx-\by_j|&\hbox{otherwise}.
\end{array}
\right.
\end{equation}
After we loop over all $S\in \mathcal{S}(\Gamma_h)$, the value $d(\bx)$ in each  $\bx\in\mathcal{N}_{\Gamma}$ is given by
\begin{equation}
d(\bx)=\min_{S_T\in \mathcal{S}(\bx)}d_T(\bx).
\end{equation}

{\it Extension phase.} During this phase, we determine both $d(\bx)$ and $u_h^e(\bx)$ for $\bx\in\mathcal{N}\setminus\mathcal{N}_\Gamma$.
To this end, the set $\mathcal{N}$ of all vertices from $\T_h$ is divided into three subsets.
A finished set $\mathcal{N}_{f}$ contains all vertices  where $d$ and $u_h^e$ have already been
defined. We initialize $\mathcal{N}_{f}=\mathcal{N}_{\Gamma}$. Initially the active set $\mathcal{N}_{a}$ contains all  the vertices, which has a neighbour in $\mathcal{N}_f$,
\begin{align*}
&\mathcal{N}_a=\{\bx\in\mathcal{N}\setminus\mathcal{N}_f~:~ \mathcal{N}(\mathcal{S}(\bx))\cap \mathcal{N}_f\neq \emptyset\},\\
&\mathcal{N}_u=\mathcal{N}\setminus(\mathcal{N}_f\cup\mathcal{N}_a).
\end{align*}
The active set is updated during the FMM and the method stops once $\mathcal{N}_a$ is empty.

For all $\bx\in\mathcal{N}_a$, the FMM iteratively computes auxiliary distance function and extension function values $\tilde{d}(\bx)$ and $\tilde{u}^e_h(\bx)$, which become final values ${d}(\bx)$ and ${u}^e_h(\bx)$ once $\bx$ leaves $\mathcal{N}_a$ and joins $\mathcal{N}_f$. The procedure is as follows. For $\bx\in\mathcal{N}_a$ we consider all $S\in \mathcal{S}(\bx)$ such that $\mathcal{N}(S)\cap\mathcal{N}_f\neq\emptyset$.
If $\mathcal{N}(S)\cap\mathcal{N}_f$ contains only one vertex $\by$, we set
\[
\tilde d_S(\bx)=d(\by)+|\bx-\by|, \quad \tilde u_{h,S}(\bx)=u_h^e(\by).
\]
If  $\mathcal{N}(S)\cap\mathcal{N}_f$ contains two or three vertices  $\{\by_j\}$, $1\leq j\leq J$, $J=2$ or $3$, then we compute
\begin{align*}
&\tilde d_S(\bx)=\left\{
\begin{array}{ll}
d(P_h \bx)+|\bx-P_h \bx|,& \hbox{if } P_h \bx\in S,\\
d(\by_{min})+|\bx-\by_{min}|,&\hbox{otherwise},
\end{array}
\right.\\
&\tilde u_{h,S}(\bx)=\left\{
\begin{array}{ll}
u_h^e(P_h \bx),& \hbox{if } P_h \bx\in S,\\
u_h^e(\by_{min}),&\hbox{otherwise},
\end{array}
\right.
\end{align*}
where $\by_{min}=\mathrm{argmin}_{1\leq j\leq m}(d(\by_j)+|\bx-\by_j|)$, and $P_h \bx$ is the orthogonal
projection of $\bx$ on the line passing through $\{\by_j\}$ (if $J=2$) or  the plane containing $\{\by_j\}$ (if $J=3$).
The value of $d(P_h \bx)$ is computed as the linear interpolation of the known values $d(\by_j)$. Now we set
\[
\begin{aligned}
\tilde{d}(\bx)&=\tilde d_{S_{min}}(\bx),\\ \tilde u_{h}(\bx)&=\tilde u_{h,S_{min}}(\bx),
\end{aligned}
\quad\text{for}~~
S_{min}=\mathrm{argmin}\{\tilde d_S(\bx)~:~S\in \mathcal{S}(\bx)~\text{and}~\mathcal{N}(S)\cap\mathcal{N}_f\neq\emptyset\}.
\]

Now we determine such vertex $\bx_{min}\in \mathcal{N}_a$  that
$
d(\bx_{min})=\min_{\bx\in\mathcal{N}_a}\tilde d(\bx).
$
and set
\[
{d}(\bx_{min})=\tilde d(\bx_{min}),\quad u^e_{h}(\bx_{min})=\tilde u_{h}(\bx_{min}).
\]
The vertex $\bx_{min}$ is now moved from the active set $\mathcal{N}_a$ to the finalized set $\mathcal{N}_f$. Based on the value of ${d}(\bx_{min})$ one checks, if any tetrahedron from  $\mathcal{S}(\bx_{min})$ may belong to $\mathcal{\widetilde{S}}(\Gamma_h)$ strip.
If such $S\in\mathcal{\widetilde{S}}(\Gamma_h)$ exists  then  $\mathcal{N}_a$ is updated by vertices  from $\mathcal{N}_u$ connected with $\bx_{\min}$. Otherwise, no new vertices are added to  $\mathcal{N}_a$. In our implementation we use the simple
criterion: If it holds
\[
d(\bx_{min})>h+L|\mathbf{w}|_{\infty}\Delta t,
\]
then we do not update $\mathcal{N}_a$  with  new vertices from $\mathcal{N}_u$.

\section{Numerical examples}\label{s_num} This section collects the results of several numerical experiments for a number of problems posed on evolving surfaces. The results demonstrate the accuracy of the trace FEM, its stability with respect to the variation of discretization parameters, and the ability to handle the case when the transport--diffusion PDE is solved on a surface undergoing topological changes.

All implementations are done in the finite element package DROPS \cite{DROPS}. The background finite element space $V_h$ consists of piecewise linear continuous finite elements. The BDF2 scheme is applied to approximate the time derivative.
At each time step, we assemble the stiffness matrix and the right-hand side by   numerical integration
over the discrete surfaces $\Gamma_h^n$. A Gaussian quadrature of degree five is used for the numerical integration on each $K\in\F_h$. The same method is used  to evaluate the finite element error. All linear algebraic  systems are solved using the  GMRES iterative method with the Gauss--Seidel preconditioner to a relative tolerance of $10^{-6}$.

The first series of experiments verifies the formal accuracy order of the method for the examples with known analytical solutions.
\medskip

\noindent{\bf Experiment 1.} We consider the transport--diffusion equation \eqref{transport}
on the unit sphere $\Gamma(t)$ moving with the constant velocity $\mathbf{w}=(0.2,0,0)$.
The initial data is given by
\[
\Gamma(0):=\{\bx\in\mathbb{R}^3~:~|\bx|=1\},\quad u|_{t=0}=1+x_1+x_2+x_3.
\]
One  easily checks that the exact solution is given by
$u(\bx,t)=1+(x_1+x_2+x_3-0.2t )\exp(-2t)$. In this and the next two experiments, we set $T=1$.

The computational domain is $\Omega=[-2,2]^3$.
We divide  $\Omega$ into  tetrahedra as follows: First we apply the uniform tessellation of $\Omega$
into cubes with side length $h$. Further the Kuhn subdivision of each cube into 6 tetrahedra is applied. This results in the shape regular background triangulation $\T_h$. The finite element level set function $\phi_h(\bx,t)$ is the nodal Lagrangian P1 interpolant for the signed distance function of $\Gamma(t)$, and
\[
\Gamma_h^n=\{\bx\in\mathbb{R}^3~:~\phi_h(\bx,t_n)=0\}.
\]
The temporal grid is uniform, $t_n=n\Delta t$. We note that in all experiments we apply the Fast Marching Method to find both distances to $\Gamma_h^n$ and $u^e$, so we never explore the explicit knowledge of the distance function for $\Gamma(t)$.

 \begin{figure}[ht!]
 \centering
  \subfigure[]{
   \includegraphics[width=0.4\textwidth]{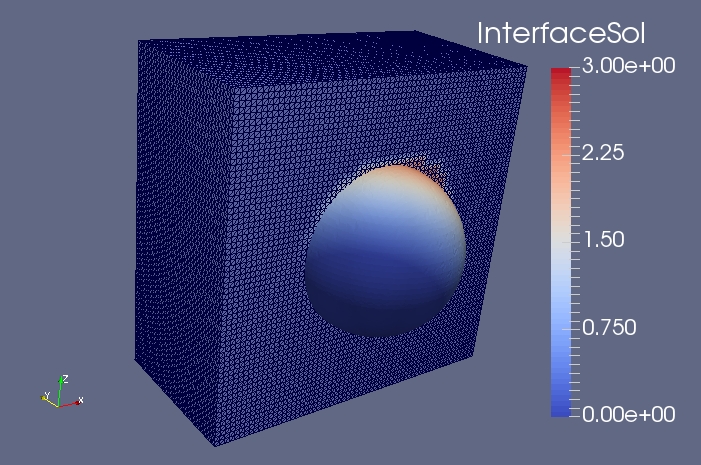}
   }
%
%
  \subfigure[]{
  \includegraphics[width=0.4\textwidth]{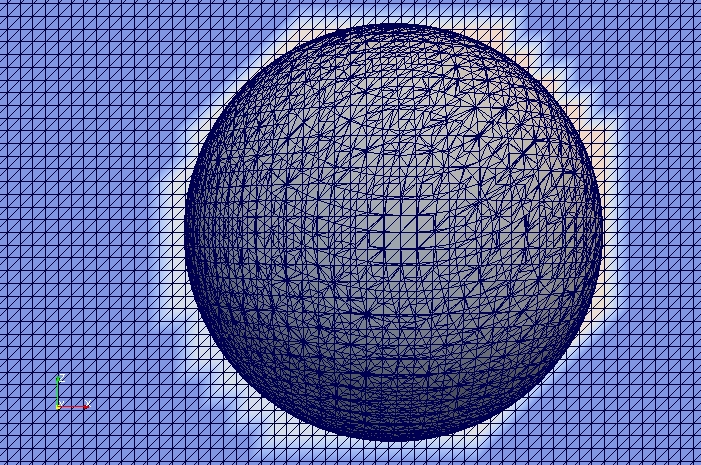}

}
 \caption{The cut of the background mesh and a part of the surface mesh for $t_n=1$. 
 Colors illustrate the solution and its extension.   
 }
 \label{fig:Exampl1}
\end{figure}

Figure~\ref{fig:Exampl1} shows the cut of the background mesh and the surface mesh colored by the computed solution at  time $t=1$.
We are interested in the $L^2(H^1)$ and $L^2(L^2)$ surface norms for the error. We compute them using the trapezoidal quadrature rule in time,
\begin{align*}
err_{L^2(H^1)} &=
\Big\{\frac{\Delta t}{2}\|\nabla_{\Gamma_h}(u^e-\pi_h u)\|_{L^2(\Gamma_h^0)}^2 +\sum_{i=1}^{N-1}\Delta t\|\nabla_{\Gamma_h}(u^e-u_h)\|_{L^2(\Gamma_h^i)}^2
\\ &\quad +\frac{\Delta t}{2}\|\nabla_{\Gamma_h}(u^e-u_h)\|_{L^2(\Gamma_h^n)}^2\Big\}^{1/2},\\
err_{L^2(L^2)} &=
\Big\{\frac{\Delta t}{2}\|(u^e-\pi_h u)\|_{L^2(\Gamma_h^0)}^2 +\sum_{i=1}^{N-1}\Delta t\|(u^e-u_h)\|_{L^2(\Gamma_h^i)}^2
\\ &\quad +\frac{\Delta t}{2}\|(u^e-u_h)\|_{L^2(\Gamma_h^N)}^2\Big\}^{1/2}.
\end{align*}
 Tables~\ref{tab:convergencerateH1} and \ref{tab:convergencerateL2} present the error norms for the Experiment~1 with various time steps $\Delta t$ and mesh sizes $h$.  If one refines both $\Delta t$ and $h$, the first order of convergence in the surface $L^2(H^1)$-norm and the second order in the surface $L^2(L^2)$-norm are clearly seen.
{\color{black} For the case of large $\Delta t$ and small $h$ the FMM extension strip $\widetilde{S}(\Gamma_h)\setminus{S}(\Gamma_h)$ becomes wider in terms of characteristic mesh size $h$, and  the accuracy of the method  diminishes. This numerical phenomenon can be noted in the top rows of Tables~\ref{tab:convergencerateH1} and \ref{tab:convergencerateL2}, where the error increases as the mesh size becomes smaller. We expect that the situation improves if one applies more accurate extension methods in   $\widetilde{S}(\Gamma_h)\setminus{S}(\Gamma_h)$. One candidate would be the normal derivative volume stabilization method from~\cite{grande2016analysis} extended to all tetrahedra in $\widetilde{S}(\Gamma_h)$.}

\begin{table}[h]\small
\caption{\small The $L^2(H^1)$-norm of the error in Experiment~1.}\label{tab:convergencerateH1}
\vspace{-0.2cm}
\begin{center}
\begin{tabular}{l|llll}
 & $h=1/2$  &$h=1/4$ &  $h=1/8$   & $h=1/16$ \\
\hline
$\Delta t=1/8$&\bf 0.96365&0.835346  &1.221340 &2.586520 \\
$\Delta t=1/16$ &0.963654&\bf 0.74794& 0.423799&0.653380\\
$\Delta t=1/32$ &0.954179&0.759253&\bf 0.37954 &0.225399\\
$\Delta t=1/64$&0.953155& 0.766650&0.381567&\bf 0.19143\\
\hline
\end{tabular}
\end{center}
\end{table}

\begin{table}[h]\small
\caption {\small The $L^2(L^2)$-norm of the error in Experiment~1.}\label{tab:convergencerateL2}
\vspace{-0.2cm}
\begin{center}
\begin{tabular}{l|llll}
 & $h=1/2$  &$h=1/4$ &  $h=1/8$   & $h=1/16$ \\
\hline
$\Delta t=1/8$&\bf 0.39351&0.192592  &0.319912 &0.691862 \\
$\Delta t=1/16$ &0.435067&\bf 0.16268& 0.057801&0.107322\\
$\Delta t=1/32$ &0.445765&0.172543&\bf 0.04013 &0.018707\\
$\Delta t=1/64$&0.448433& 0.175145&0.041875&\bf 0.01040\\
\hline
\end{tabular}
\end{center}
\end{table}

\begin{table}[h]\small
{\color{black}\caption{\small Averaged CPU times per each time step of the method in Experiment~1}\label{tab:computeTime}
\vspace{-0.2cm}
\begin{center}
\begin{tabular}{c|cclll}
& Active d.o.f.  & Extra d.o.f.     & $T_{assemb}$ & $T_{solve}$ & $T_{ext}$ \\
\hline
$h_0=1/2,~\Delta t_0=1/8$      &31 &8     &0.0038&0.0004&0.0012  \\
$h=h_0/2,~\Delta t=\Delta t_0/2$& 104 &24 &0.0160&0.0021&0.0041 \\
$h=h_0/4,~\Delta t=\Delta t_0/4$  &452&63 &0.0708&0.0087&0.0195\\
$h=h_0/8,~\Delta t=\Delta t_0/8$&1880&170 &0.3814&0.0374&0.0906\\
\hline
\end{tabular}
\end{center}
}
\end{table}
{\color{black}
Table~\ref{tab:computeTime} shows the breakdown of the computational costs  of the method into 
the averaged CPU times for assembling stiffness matrices, solving resulting linear algebraic systems,
and performing the extension to $\widetilde{S}(\Gamma_h)\setminus{S}(\Gamma_h)$ by FMM.  Since 
the surface evolves, all the statistics slightly vary in time, and so the table shows averaged numbers per one time step. ``Active d.o.f.'' is the dimension of the linear algebraic system, i.e. the number of bulk finite element nodal values tailored to tetrahedra from $S(\Gamma_h)$. ``Extra d.o.f.'' is the
number of mesh nodes in  $\widetilde{S}(\Gamma_h)\setminus\overline{{S}(\Gamma_h)}$, these are all nodes where extension is computed by  FMM. The averaged CPU times demonstrate  optimal  or close to the optimal scaling with respect to the number of degrees of freedom. As common for a finite element method, the most time consuming part is the assembling of the stiffness matrices. The costs of  FMM are modest compared to the assembling time, and $T_{solve}$ indicates that using a preconditioned Krylov subspace method is the efficient approach to solve linear algebraic systems (no extra stabilizing terms were added to
the FE formulation for improving its algebraic properties).
}


\medskip

\noindent{\bf Experiment~2.} 
The setup of this experiment is similar to the previous one.  The transport velocity is given by
$\mathbf{w}=(-2\pi x_2, 2\pi x_1,0)$. Initially, the sphere is set off the center of the domain:
The initial data is given by
\[
\Gamma(0):=\{\bx\in\mathbb{R}^3~:~|\bx-\bx_0|=1\},\quad u|_{t=0}=1+(x_1-0.5)+x_2+x_3,
\]
with $\bx_0=(0.5,0,0)$. Now $\bw$ revolves the sphere around the center of domain without changing its shape.
One  checks that the exact solution to \eqref{transport} is given by
\[
u(\bx,t)=(x_1(\cos(2\pi t)-\sin(2\pi t))+x_2(\cos(2\pi t)+\sin(2\pi t))+x_3+0.5 )\exp(-2t).
\]

{
\begin{table}\small
\caption{\small The $L^2(H^1)$-norm of the error in Experiment~2.}\label{tab:ErrorH1_Ex2}
\vspace{-0.2cm}
\begin{center}
\begin{tabular}{l|llll}
 & $h=1/2$  &$h=1/4$ &  $h=1/8$   & $h=1/16$ \\
\hline
{\color{black}$\Delta t=1/32$}& 0.978459 &1.931081&3.740840&4.048480\\
$\Delta t=1/64$&\bf 0.90425& 0.690963& 0.813820& 1.066030 \\
$\Delta t=1/128$ &0.901234 &\bf 0.64014 &0.348516 & 0.300654\\
$\Delta t=1/256$ &0.901443& 0.640055& \bf 0.32352& 0.171101\\
$\Delta t=1/512$&0.901631& 0.641018& 0.323199& \bf 0.16286\\
\hline
\end{tabular}
\end{center}
\end{table}
\begin{table}\small
\caption{\small The $L^2(L^2)$-norm of the error in Experiment~2.}\label{tab:ErrorL2_Ex2}
\vspace{-0.2cm}
\begin{center}
\begin{tabular}{l|llll}
 & $h=1/2$  &$h=1/4$ &  $h=1/8$   & $h=1/16$  \\
\hline
{\color{black}$\Delta t=1/32$}& 0.294122 &0.548567&0.975485&0.958447\\
$\Delta t=1/64$&\bf 0.27244& 0.120061& 0.152520& 0.175920 \\
$\Delta t=1/128$ &0.279106&\bf 0.10451& 0.034962& 0.037085\\
$\Delta t=1/256$ & 0.279744 &0.105975& \bf 0.02699& 0.010692\\
$\Delta t=1/512$&0.279811& 0.106116& 0.026444&\bf 0.00736\\
\hline
\end{tabular}
\end{center}
\end{table}
}

 Tables~\ref{tab:ErrorH1_Ex2} and \ref{tab:ErrorL2_Ex2} show the error norms for the Experiment~2 with various time steps $\Delta t$ and mesh sizes $h$. If one refines both $\Delta t$ and $h$, the first order of convergence in the surface $L^2(H^1)$-norm and the second order in the surface $L^2(L^2)$-norm are again observed.
{\color{black}Note that the transport velocity $\|\mathbf{w}\|_{\infty}\approx 9.42$ in this experiment scales differently compared to Experiment~1. Therefore, we consider smaller $\Delta t$ to obtain meaningful results.}
\bigskip

\noindent{\bf Experiment~3.} In this experiment, we consider a shrinking sphere and solve \eqref{transport} with a source term on the right-hand side.
The bulk velocity field is given by
$
\mathbf{w}={ -{\frac12 e^{-t/2}}}\mathbf{n},
$
where $\mathbf{n}$ is the unit outward normal on $\Gamma(t)$. $\Gamma(0)$ is the unit sphere. One computes ${\Div}_\Gamma \bw = -1$.
The prescribed analytical solution  $u(\bx,t)=(1+x_1x_2x_3) e^{t}$ solves \eqref{transport} with the  right-hand side
$
f(\bx,t)=(-1.5 e^{t}+{12}e^{2t})x_1x_2x_3.
$

{
\begin{table}\small
\caption{\small The $L^2(H^1)$-norm of the error in Experiment~3.}\label{tab:ErrorH1_Ex3}\vspace{-0.2cm}
\begin{center}
\begin{tabular}{l|llll}
 & $h=1/4$  &$h=1/8$ &  $h=1/16$   & $h=1/32$ \\
\hline
$\Delta t=1/16$&\bf 0.48893 &0.311146 & 0.170104 &0.088521 \\
$\Delta t=1/32$ &0.481896& \bf  0.30859 & 0.168635& 0.087013\\
$\Delta t=1/64$ &0.478675& 0.307416 &\bf 0.16801& 0.086513\\
$\Delta t=1/128$&0.477226& 0.306872& 0.167747&\bf 0.08634\\
\hline
\end{tabular}
\end{center}
\end{table}

\begin{table}\small
\caption{\small  The $L^2(L^2)$-norm of the error in Experiment~3.}\label{tab:ErrorL2_Ex3}\vspace{-0.2cm}
\begin{center}
\begin{tabular}{l|llll}
 & $h=1/4$  & $h=1/8$ &  $h=1/16$   & $h=1/32$ \\
\hline
$\Delta t=1/16$&\bf 0.12237& 0.0468812& 0.0224705& 0.0178009 \\
$\Delta t=1/32$ &0.116763& \bf 0.040745& 0.0130912& 0.0060213 \\
$\Delta t=1/64$ & 0.115589& 0.0396511&\bf 0.011517& 0.0035199\\
$\Delta t=1/128$&0.115336& 0.0394023& 0.0112094& \bf 0.003038\\
\hline
\end{tabular}
\end{center}
\end{table}
}

 Tables~\ref{tab:ErrorH1_Ex3} and \ref{tab:ErrorL2_Ex3} show the error norms for various time steps $\Delta t$ and mesh sizes $h$. If one refines both $\Delta t$ and $h$, the first order of convergence in the surface $L^2(H^1)$-norm and the second order in the surface $L^2(L^2)$-norm is observed for the example of the shrinking sphere.

\bigskip

\noindent{\bf Experiment~4.} In this example, we consider a surface  transport--diffusion problem as in \eqref{transport}  on a more complex moving manifold. The initial manifold and concentration are given ({as in \cite{Dziuk88}}) by
$
\Gamma(0)=\{ \, \bx\in \mathbb{R}^3~:~ (x_1 -x_3^2)^2+x_2^2+x_3^2=1\, \}, \quad u_0(\bx)=1+x_1 x_2 x_3.
$
The velocity field that transports the surface is
$$\mathbf{w}(\bx,t)=\big(0.1x_1 \cos(t),0.2x_2 \sin(t),0.2x_3\cos(t)\big)^T.$$

\begin{figure}[ht!]
 \centering
  \subfigure[$t=0$]{
   \includegraphics[width=1.8in]{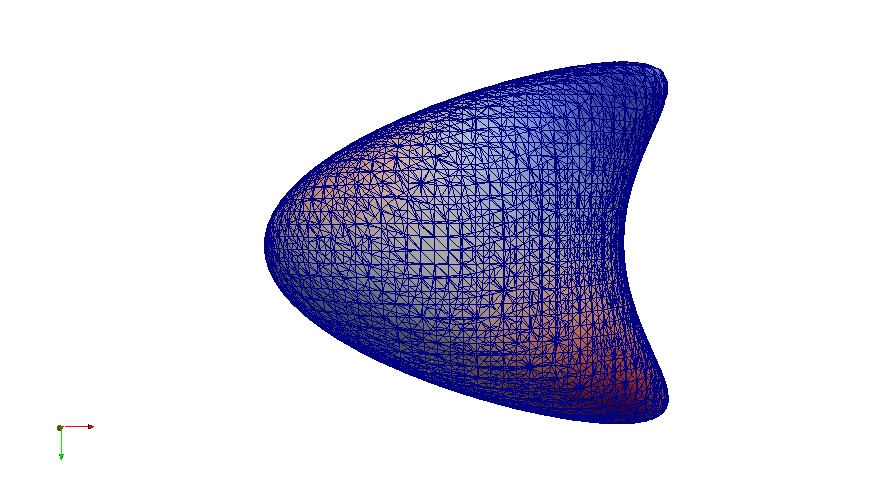}

   }
     \subfigure[$t=0.1$]{
   \includegraphics[width=1.8in]{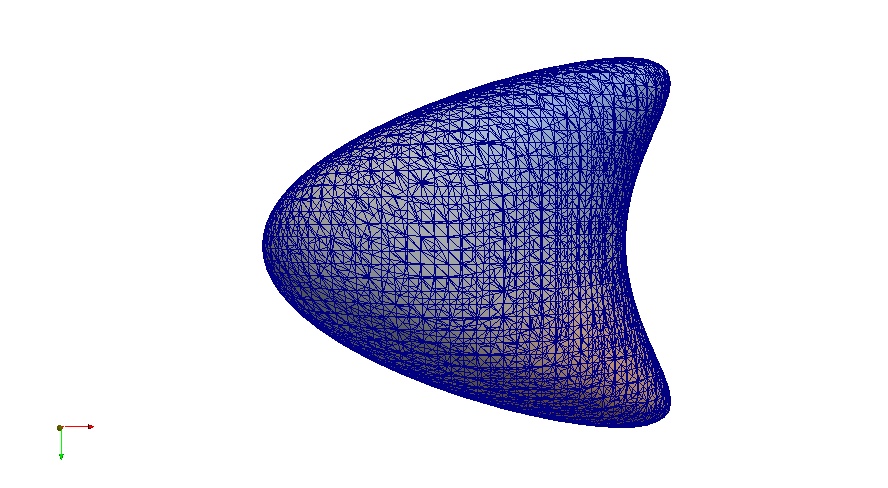}
   }
  \subfigure[$t=1$]{
  \includegraphics[width=1.8in]{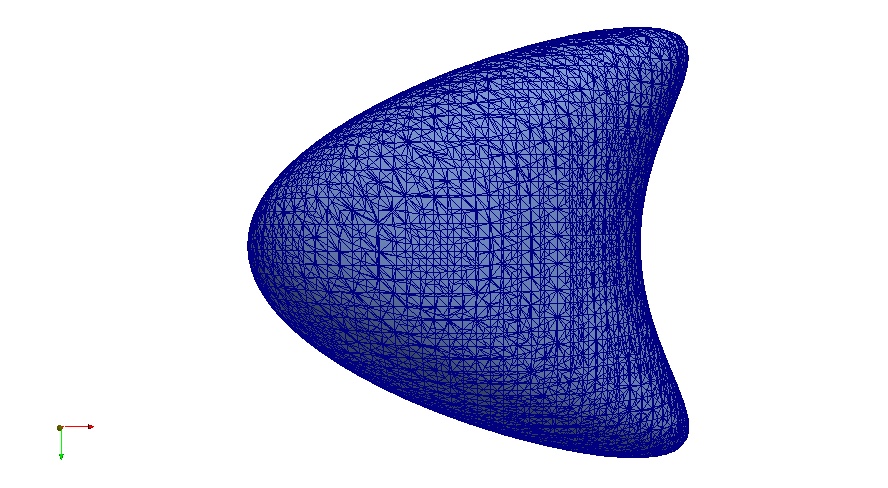}

}
  \subfigure[$t=2$]{
  \includegraphics[width=1.8in]{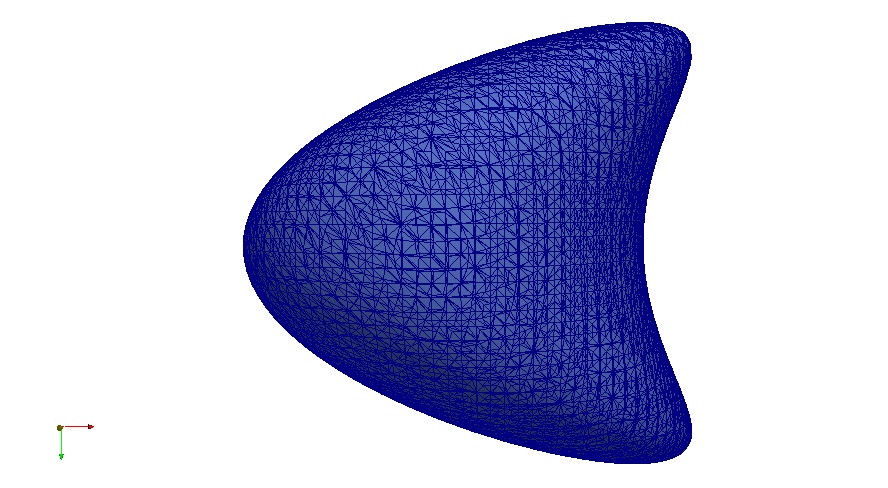}

}
  \subfigure[$t=4$]{
  \includegraphics[width=1.8in]{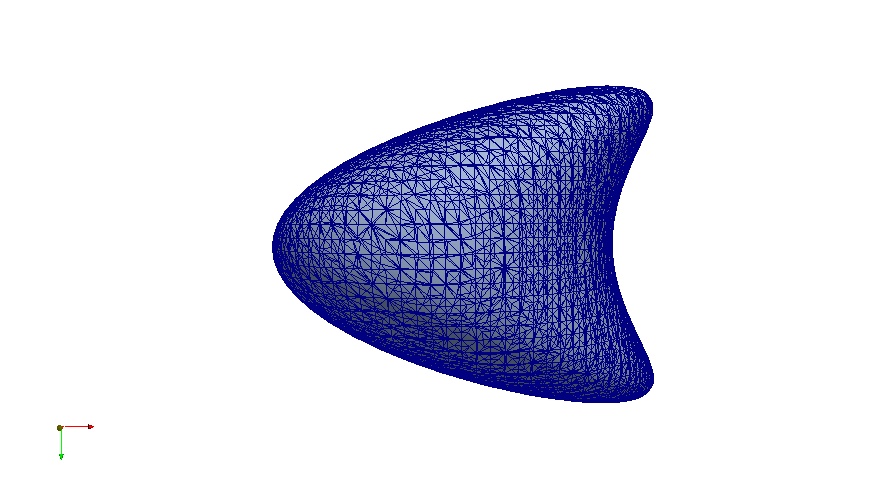}

}
  \subfigure[$t=6$]{
  \includegraphics[width=1.8in]{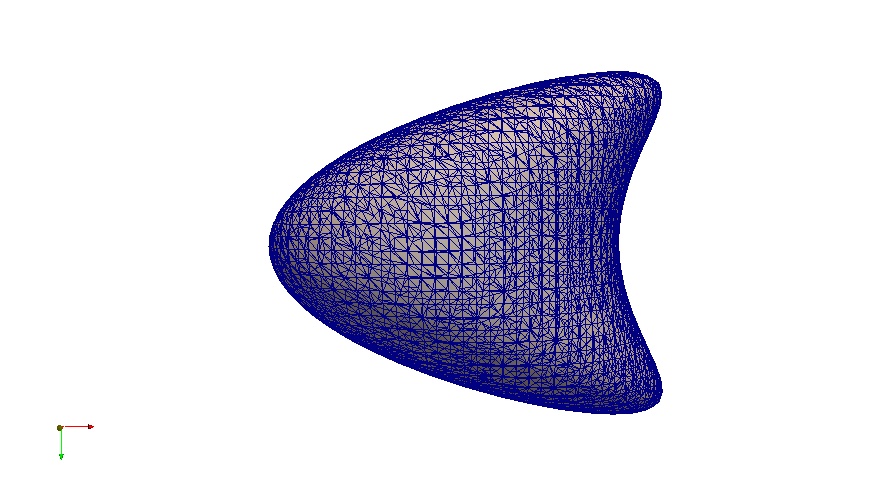}

}
 \caption{Snapshots  of the surface, surface mesh and the computed solution from Experiment~4.
 }
 \label{fig:Exampl4}
\end{figure}

\begin{figure}[ht!]
 \centering
  \subfigure[]{
   \includegraphics[width=2.8in]{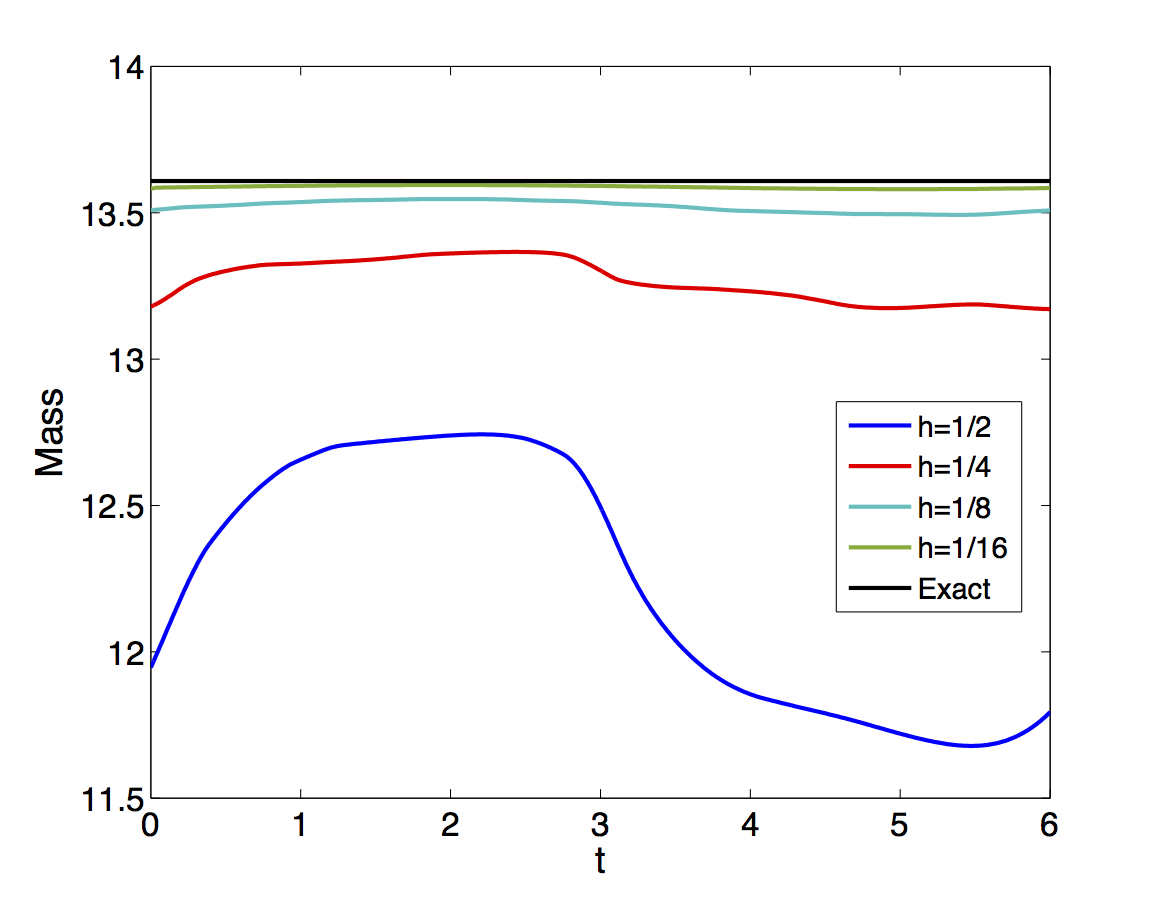}
   }
  \subfigure[]{
  \includegraphics[width=2.8in]{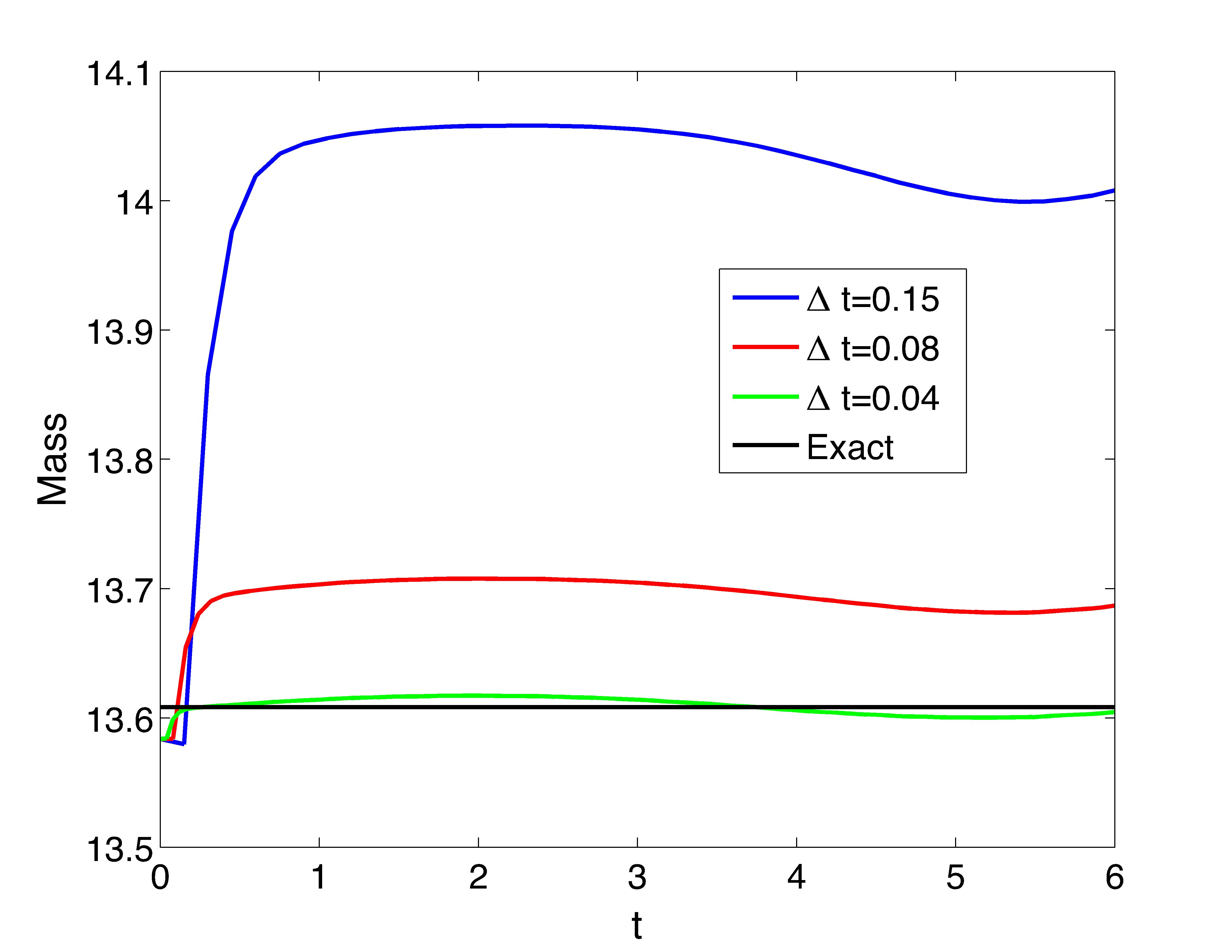}
}
 \caption{Total mass evolution for the finite element solution in Experiment 4.}
 \label{fig:mass}
\end{figure}

{\color{black}We  compute the problem until $T=6$}. In this example,
the total mass
$
M(t)=\int_{\Gamma(t)}u(\cdot,t) \, \ds
$
is conserved  and equal to
$M(0)=|\Gamma(0)|\approx 13.6083.$
 We check how well the discrete quantity $
M_h(t)=\int_{\Gamma_h(t)}u_h(\cdot, t) \, \ds
$ is conserved.
In Figure~\ref{fig:mass} (left) we plot $M_h(t)$  for different mesh sizes $h$ and a fixed time step $\Delta t=0.01$. The error in the total mass at $t=T$ is equal to  1.8142, 0.4375,0.1006 and 0.0239 (for mesh sizes as in Figure~\ref{fig:mass} (left)).
In Figure~\ref{fig:mass} (right) we plot $M_h(t)$ for different time steps $\Delta t$ and a fixed mesh size $h=1/16$.  The error in the total mass at $t=T$ is equal  0.3996, 0.0785, 0.0038 (for time steps as in Figure~\ref{fig:mass} (right)).
The error in the mass conservation is consistent with the expected second order accuracy in time and space.

If one is interested in the exact mass conservation on the discrete level, then one may enforce
$M_h(t_n)=M_h(0)$ as a side constrain in the finite element formulation \eqref{transportFDFE} with the help of the scalar Lagrange multiplier, see \cite{hansbo2016cut}. Here we used the error reduction in total mass as an \textit{indicator} of the method convergence order for the case, when the exact solution is not available.
\medskip

\noindent{\bf Experiment~5.} In this test problem from~\cite{GOReccomas}, one solves the transport--diffusion equation \eqref{transport} on an evolving surface $\Gamma(t)$ which undergoes a change of topology and experiences a local singularity.
The computational domain is $\Omega=(-3,3)\times(-2,2)^2$, $t \in [0,1]$.   The evolving surface is the zero level of the level set function $\phi$ defined as:
\[
  \phi(\bx,t) = 1 - \frac{1}{\| \bx -c_+(t)\|^3} - \frac{1}{\| \bx -c_-(t)\|^3},
\]
with  $c_\pm(t)= \pm\frac32(t - 1, 0 , 0)^T$, $t \in [0,1]$.
For $t=0$ and $\bx \in B(c_+(0);1)$ one has $\| \bx -c_+(0)\|^{-3} =1$ and $\|\bx -c_-(0)\|^{-3} \ll 1$. For $t=0$ and  $ \bx \in B(c_-(0);1)$ one has $\| \bx -c_+(0)\|^{-3}  \ll 1$ and $\| \bx -c_-(0)\|^{-3} =1$. Hence,  the initial configuration $\Gamma(0)$ is  close to two balls of radius $1$, centered at $\pm(1.5, 0, 0)^T$. For $t=1$ the surface $\Gamma(1)$ is the ball around $0$ with radius $2^{1/3}$. For $t >0 $ the two spheres approach each other until time $\tilde t= 1-\tfrac23 2^{1/3}\approx 0.160$, when they touch at the origin. For $t \in (\tilde t,1]$ the surface $\Gamma(t)$ is simply connected and deforms into the sphere $\Gamma(1)$.

\begin{figure}[ht!]
 \centering
  \subfigure[$t=0$]{
   \includegraphics[width=1.8in]{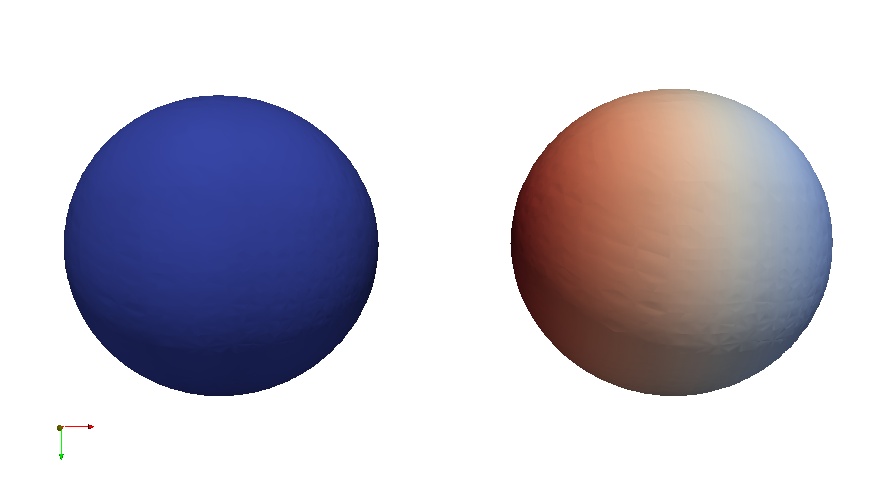}
   }
  \subfigure[$t=1/8$]{
  \includegraphics[width=1.8in]{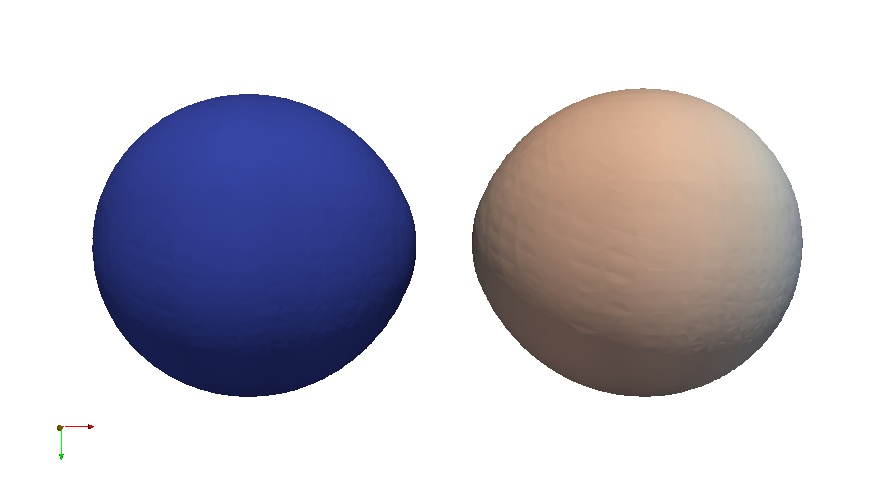}
}
  \subfigure[$t=1/4$]{
   \includegraphics[width=1.8in]{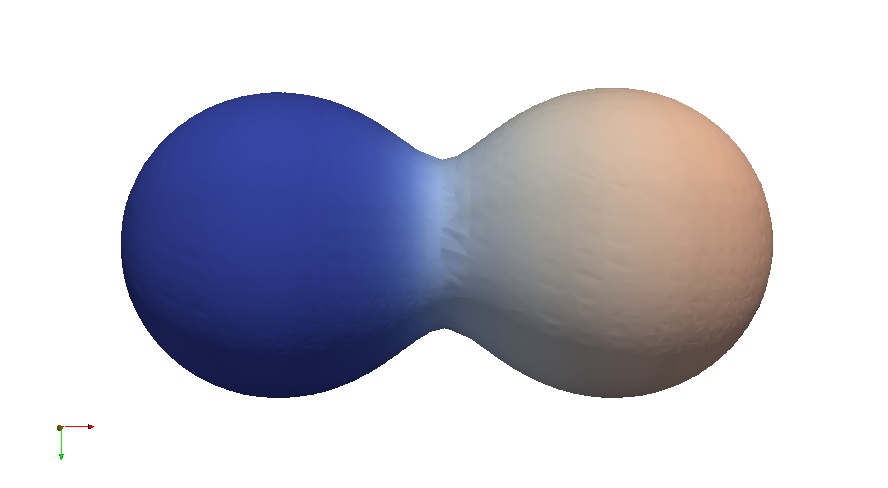}
   }
  \subfigure[$t=3/8$]{
  \includegraphics[width=1.8in]{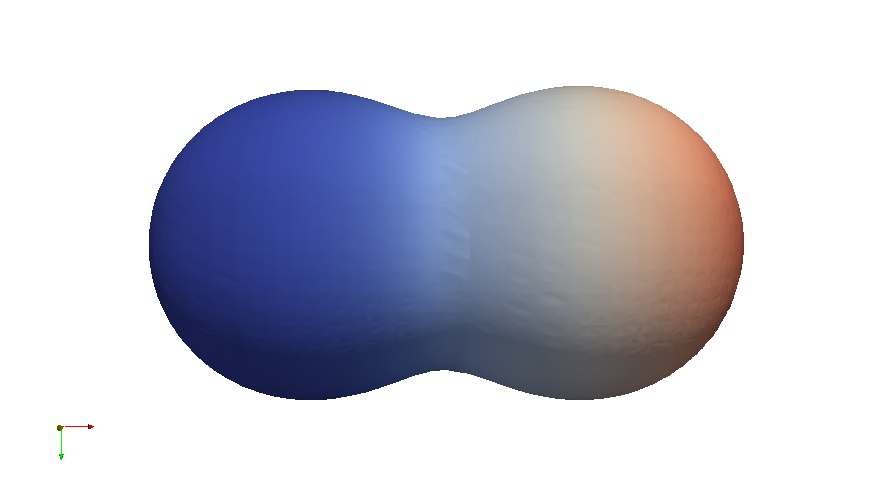}
}
  \subfigure[$t=1/2$]{
  \includegraphics[width=1.8in]{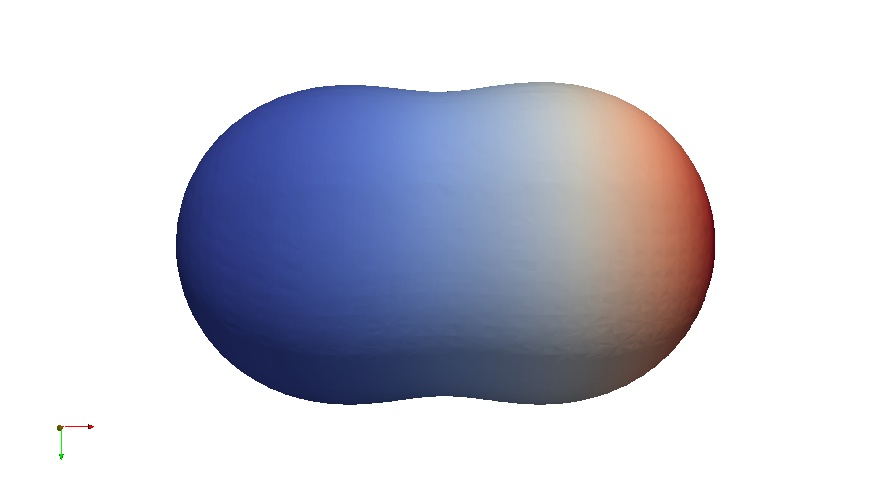}
}
  \subfigure[$t=5/8$]{
  \includegraphics[width=1.8in]{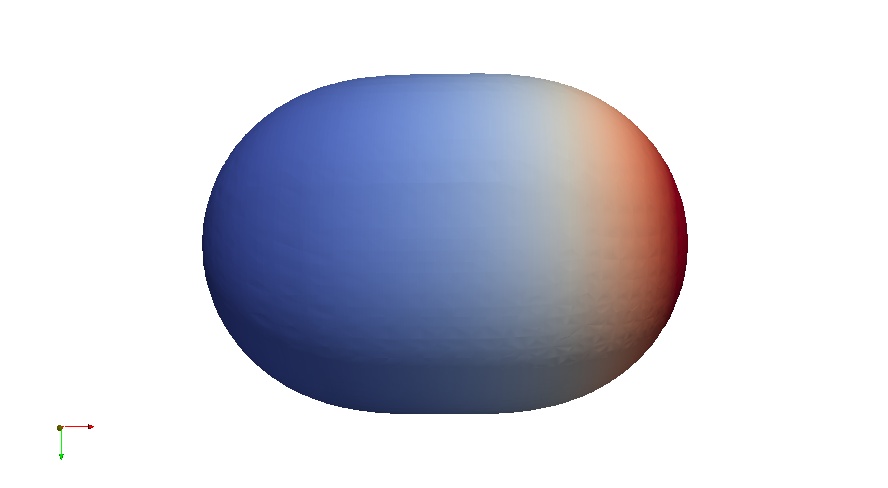}
}
  \subfigure[$t=3/4$]{
  \includegraphics[width=1.8in]{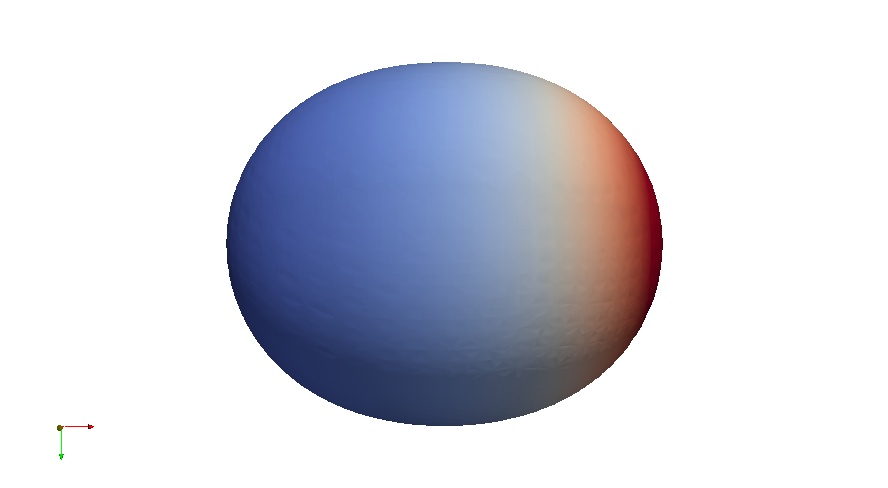}
}
  \subfigure[$t=7/8$]{
  \includegraphics[width=1.8in]{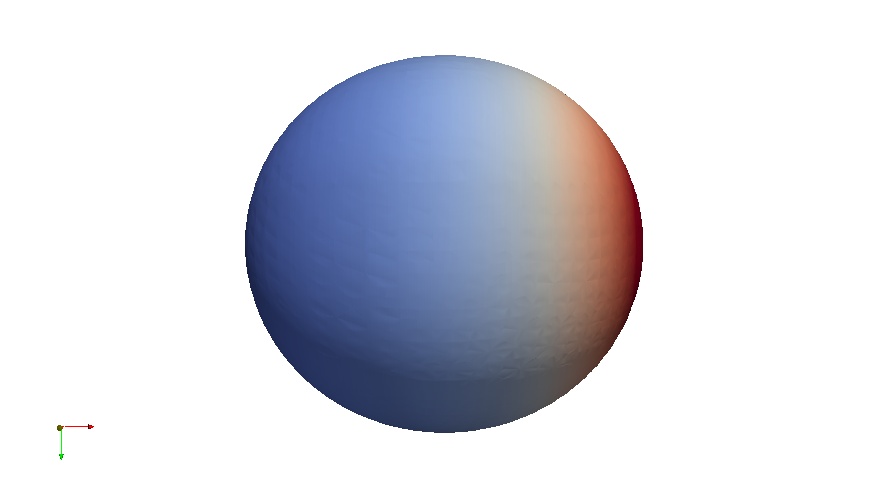}
}
  \subfigure[$t=1$]{
  \includegraphics[width=1.8in]{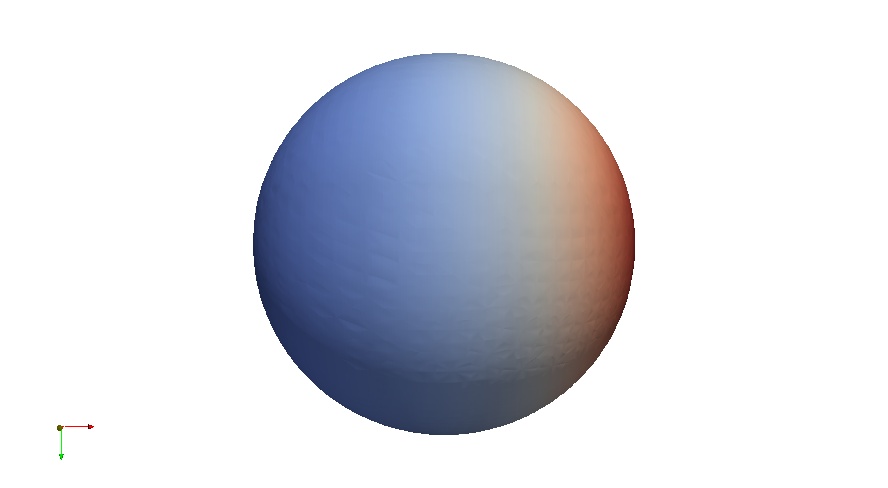}
}
 \caption{
 Snapshots of discrete solution in Experiment~5 with $h=1/16$, $\Delta t=1/128$.}
 \label{fig:solutions}
\end{figure}

In the vicinity of $\Gamma(t)$, the gradient $\nabla\phi$ and the time derivative $\partial_t\phi$ are well-defined and given by simple algebraic expressions. The normal velocity field, which transports $\Gamma (t)$, can be computed (cf. \cite{GOReccomas}) to be
\[
  \bw= -\frac{\partial_t\phi}{\lvert\nabla\phi\rvert^2}\nabla\phi.
\]
The initial value of $u$ is given by
\[
 u_0(\bx)=\left\{\begin{array}{ll}
3-x_1 & \hbox{for }x_1\geq 0;\\
0&\hbox{otherwise.}
 \end{array}
 \right.
\]

\begin{figure}[ht!]
 \centering
  \subfigure[$t=0.1484$]{
   \includegraphics[width=2.1in]{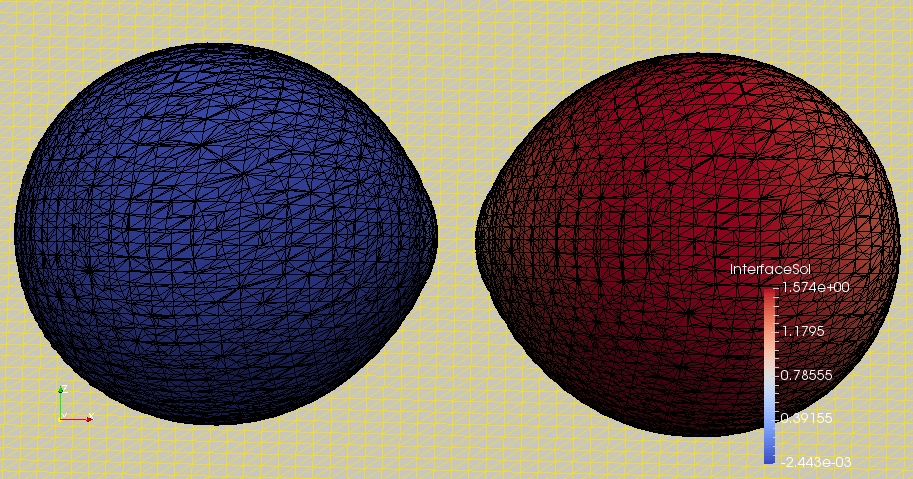}
   }
  \subfigure[$t=0.1484$ (zoomed-in)]{
  \includegraphics[width=2.1in]{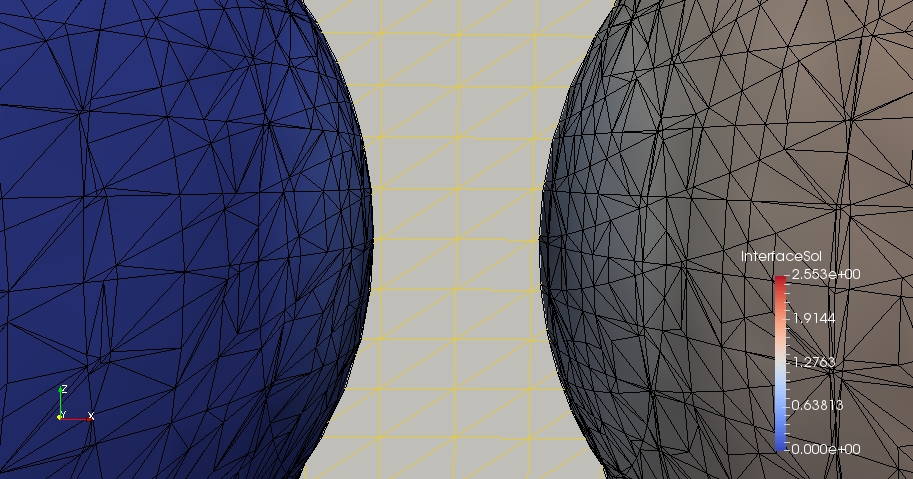}
}
  \subfigure[$t=0.15625$ ]{
   \includegraphics[width=2.1in]{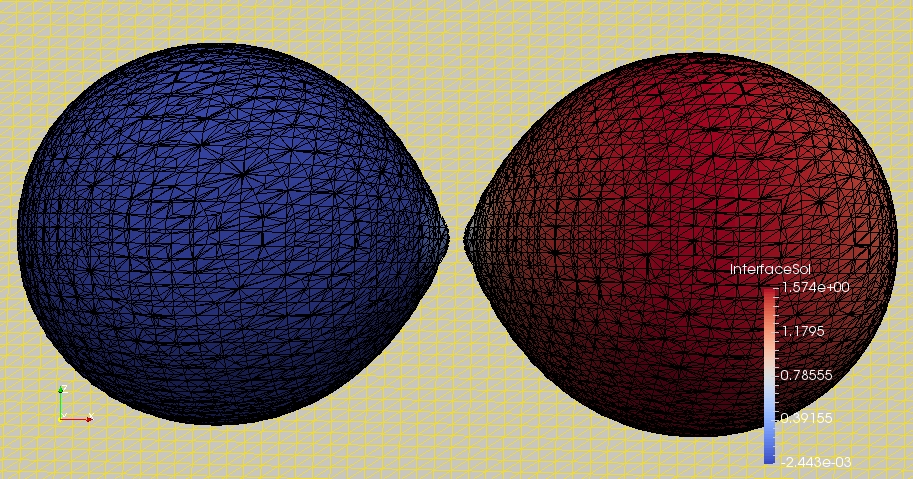}
   }
  \subfigure[$t=0.15625$ (zoomed-in)]{
  \includegraphics[width=2.1in]{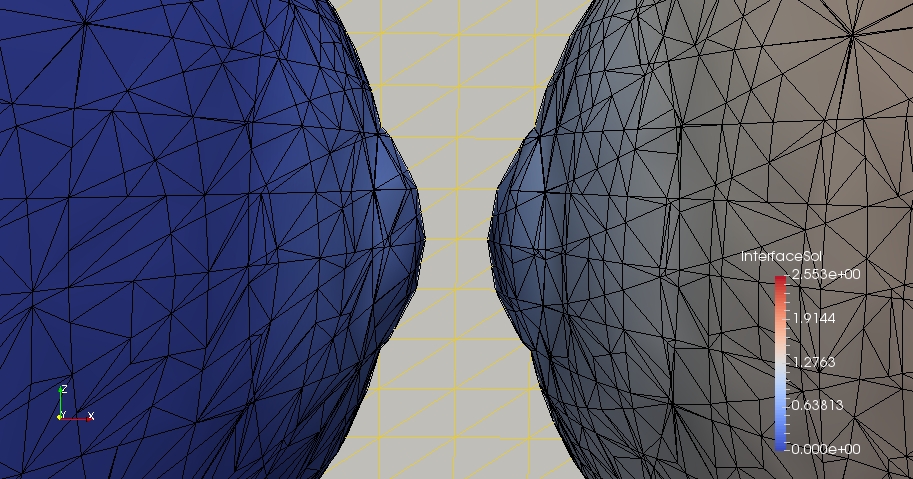}
}
  \subfigure[$t=0.1641$]{
  \includegraphics[width=2.1in]{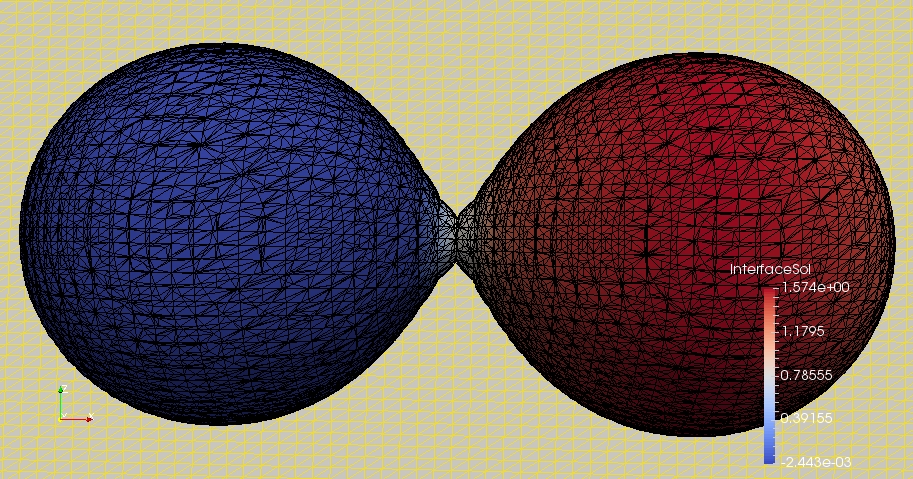}
}
  \subfigure[$t=0.1641$ (zoomed-in)]{
  \includegraphics[width=2.1in]{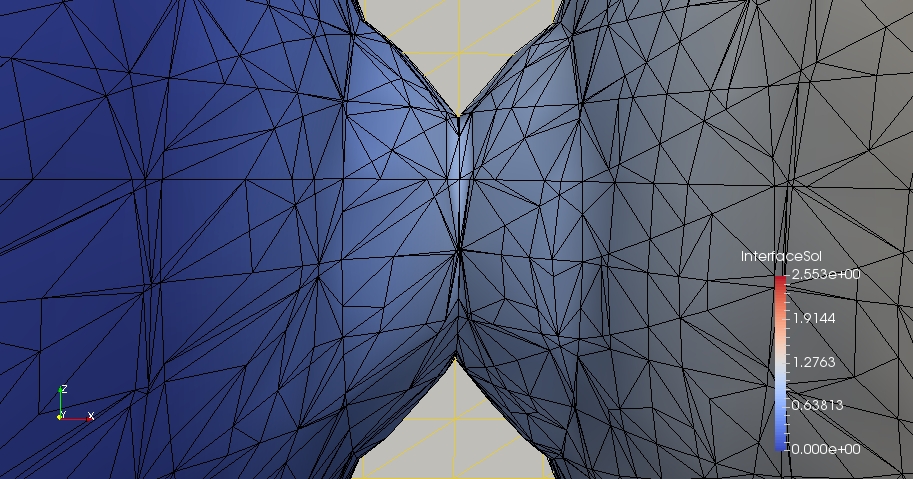}
}
  \subfigure[$t=0.171875$]{
  \includegraphics[width=2.1in]{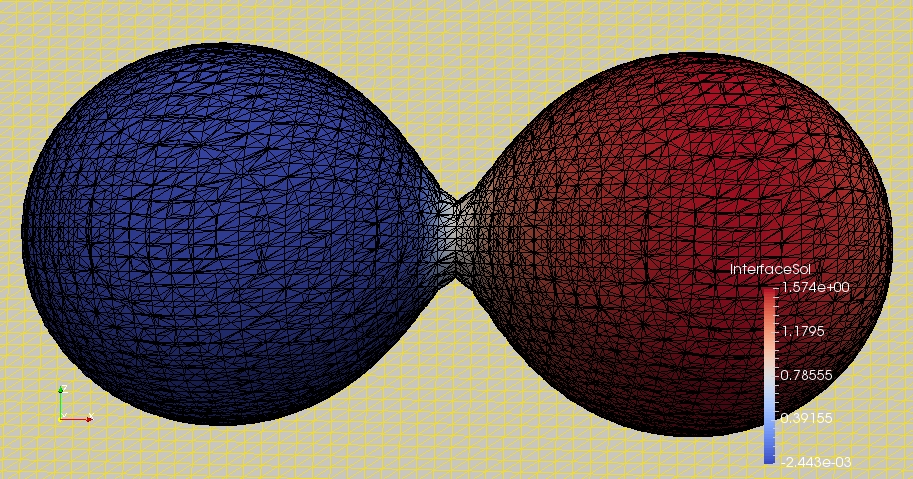}
}
  \subfigure[$t=0.171875$ (zoomed-in)]{
  \includegraphics[width=2.1in]{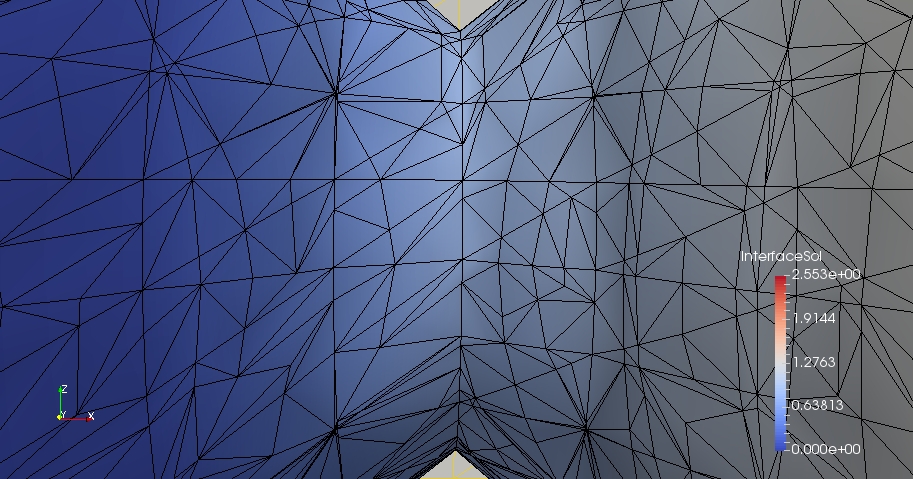}
}
 \caption{The computed solution and surface meshes close to the time of collision Example~5.}
 \label{fig:surfaceMeshMerge}
\end{figure}

 \begin{figure}[ht!]
 \centering
 \subfigure[]{
  \includegraphics[width=2.3in]{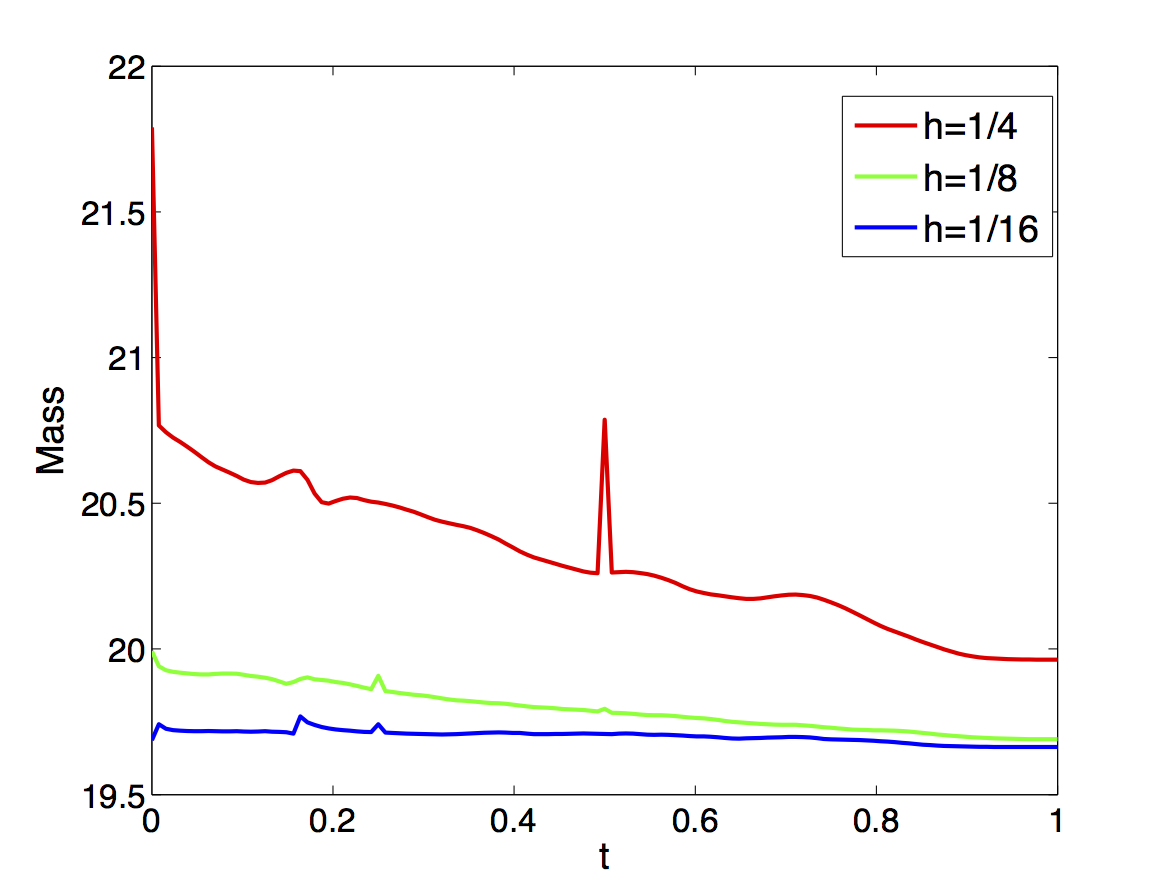}
  }
   \subfigure[]{
  \includegraphics[width=2.3in]{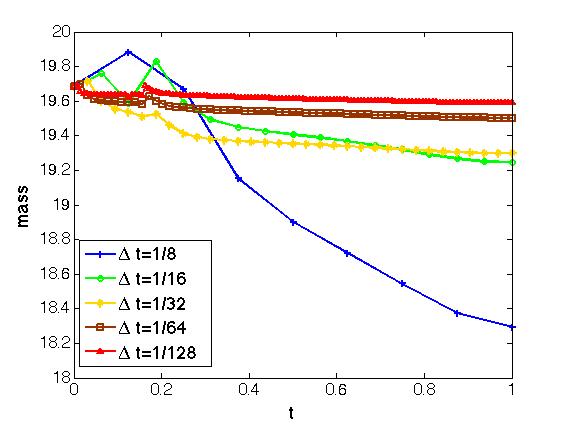}
  }
 \caption{Total mass evolution for the finite element solution in Experiment~5.}
 \label{fig:mass_merge}
\end{figure}

In Figure~\ref{fig:solutions} we show a few snapshots of the surface and the computed surface solution on for the background tetrahedral mesh with $h=1/16$  and $\Delta t =1/128$. The surfaces $\Gamma_h^n$ close to the time of collision are illustrated in Figure~\ref{fig:surfaceMeshMerge}. The suggested variant of the trace FEM handles the geometrical singularity without any difficulty.
{\color{black} It is clear that the computed extension $u^e$ in this experiment is not smooth in a neighborhood of the singularity, and so formal analysis of the consistency of the method is not directly applicable to this case. However, the closest point extension is well defined and numerical results suggest that this is sufficient for the method to be stable.}
Similar to the previous example, we compute the total discrete mass $M_h(t)$ on $\Gamma_h^n$. This can be used as a measure of accuracy. The evolution of $M_h(t)$ for varying $h$ and $\Delta t$ is shown in  Figure~\ref{fig:mass_merge}. The convergence of the quantity is obvious. Finally, we note that in this experiment we observed the stable numerical behaviour of the method
for any combinations of the mesh size and time step we tested, including $\Delta t = \frac18$, $h=\frac1{16}$ and $\Delta t = \frac1{128}$, $h=\frac1{4}$.

\section{Conclusions}
We studied a new fully Eulerian unfitted finite element method for solving PDEs posed on evolving surfaces. The
method combines three computational techniques in a modular way: a finite difference approximation in time,
a finite element method on stationary surface, and an extension of finite element functions from a surface to a neighborhood of the surface. All of these three computational  techniques have been intensively studied in the literature,
and so well established methods can be used. In this paper, we used the trace piecewise linear finite element method -- the higher order variants of this method are also available in the literature~{\color{black} \cite{grande2016analysis}} -- for the spacial discretization and a variant of the fast marching method to construct suitable extension. We observed stable and second order accurate numerical results in a several experiments including  one experiment with two colliding spheres. The rigorous  error analysis of the method is a subject of further research.
\bibliography{literatur}{}
\bibliographystyle{abbrv}
\end{document}